\documentclass[10pt]{article}
\usepackage{amsmath,amssymb,amsthm}
\usepackage{color}
\usepackage[nohead,margin=1.0in]{geometry}
\usepackage{booktabs}
\usepackage{algorithm,algorithmic}
\usepackage{verbatim}
\usepackage{url}
\usepackage{graphicx}
\usepackage{epstopdf}
\graphicspath{{./pics/}}
\allowdisplaybreaks

 \usepackage{array}
\theoremstyle{plain}
\newtheorem{theorem}{Theorem}[section]
\newtheorem{remark}{Remark}[section]
\newtheorem{definition}{Definition}[section]
\newtheorem{lemma}{Lemma}[section]
\newtheorem{assumption}{Assumption}[section]
\newtheorem{proposition}{Proposition}[section]
\newtheorem{corollary}{Corollary}[section]
\newtheorem{example}{Example}[section]
\numberwithin{equation}{section}

\renewcommand{\d}{\mathrm{d}}

\title{Reconstruction of a Space-Time Dependent Source in Subdiffusion Models via a Perturbation Approach\thanks{The work of B.J. is partially supported by UK EPSRC grant EP/T000864/1, that of Y.K. by the French National Research Agency ANR
(project MultiOnde) grant ANR-17-CE40-0029, and that of Z.Z. by Hong Kong RGC grant (No. 15304420).}}
\author{Bangti Jin\thanks{Department of Computer Science, University College London, Gower Street, London WC1E 6BT, UK (\texttt{b.jin@ucl.ac.uk})}
\and Yavar Kian\thanks{Centre de Physique Th\'{e}orique  (CPT),  UMR-7332, Aix Marseille Universit\'{e}, Campus de Luminy, Case 907,
13288 Marseille cedex 9, France (yavar.kian@univ-amu.fr)} \and Zhi Zhou\thanks{Department of Applied Mathematics, The Hong Kong Polytechnic University, Hung Hom, Kowloon, Hong Kong (zhi.zhou@polyu.edu.hk)}}

\begin{document}

\maketitle
\begin{abstract}
In this article we study two inverse problems of recovering a space-time dependent source component from
the lateral boundary observation in a subidffusion model. The mathematical model involves a
Djrbashian-Caputo fractional derivative of order $\alpha\in(0,1)$ in time, and a second-order elliptic
operator with time-dependent coefficients. We establish a well-posedness and a conditional stability
result for the inverse problems using a novel perturbation argument and refined regularity estimates
of the associated direct problem. Further, we present a numerical algorithm for efficiently and
accurately reconstructing the source component, and provide several two-dimensional numerical results
showing the feasibility of the recovery.
\end{abstract}

\section{Introduction}

This work is concerned with inverse source problems (ISPs) of identifying a space-time dependent component of the source in
the subdiffusion model in a cylindrical domain from the lateral Cauchy data on a part of the boundary. Let $d\geq 2$, $\Omega=\omega\times(-\ell,\ell)$,
$\omega\subset \mathbb{R}^{d-1}$ be an open bounded domain with a $C^2$ boundary, and fix $T>0$ the final time. For any $x\in\Omega$,
we write $x=(x',x_d)$, with $x'\in \omega$ and $x_d\in(-\ell,\ell)$. For $m=0,1$, we consider the following initial boundary value problem
for the function $u$:
\begin{equation}\label{eqn:fde}
  \left\{\begin{aligned}
    \partial_t^\alpha u + \mathcal{A}(t) u &= F,\quad \mbox{in }\Omega\times (0,T),\\
    u(x,0) & = 0,\quad \mbox{in }\Omega,\\
    \partial_{x_d}^mu(x',\ell,t) & = 0, \quad \mbox{on }\omega\times(0,T),\\
    \partial_{x_d} u(x',-\ell,t) & = 0, \quad \mbox{on }\omega\times(0,T),\\
   u(x,t) & = 0,\quad \mbox{on } \partial\omega\times(-\ell,\ell) \times (0,T).
  \end{aligned}  \right.
\end{equation}
In the model \eqref{eqn:fde}, the order $\alpha\in (0,1)$ is fixed, and the notation $\partial_t^\alpha u$ denotes the so-called Djrbashian-Caputo
fractional derivative of order $\alpha$ in time, which, for $\alpha\in(0,1)$, is defined by \cite[p. 92]{KilbasSrivastavaTrujillo:2006}
\begin{equation*}
  \partial_t^\alpha u(t) = \frac{1}{\Gamma(1-\alpha)}\int_0^t(t-s)^{-\alpha}u'(s)\d s,
\end{equation*}
where $\Gamma(z)=\int_0^\infty s^{z-1}e^{-s}\d s$ for $\Re(z)>0$ denotes Euler's Gamma function (the notation
$\Re$ denotes taking the real part of a complex number $z$). When the order $\alpha$ approaches
$1^-$, the fractional derivative $\partial_t^\alpha u$ recovers the usual first-order derivative $u'(t)$,
and accordingly, the model coincides with the standard diffusion equation. $\mathcal{A}(t)$ is a
time-dependent second-order strongly elliptic operator, defined by
\begin{equation*}
  \mathcal{A}(t) u (x) = -\sum_{i,j=1}^{d} \partial_{x_i}(a_{ij}(x,t)\partial_{x_j}u) +q(t)u,\quad x\in\Omega,
\end{equation*}
where $a=[a_{ij}]_{i,j=1}^d\in C(\overline{\Omega}\times[0,T];\mathbb{R}^{d\times d})$ is a symmetric matrix-valued function
and satisfies suitable regularity conditions given in Assumption \ref{ass:a} below, and $q\in C^1([0,T];L^{\infty}(\Omega))$
is nonnegative.

The model \eqref{eqn:fde} has received much attention in recent years, known under the name of ``subdiffusion'' or ``time-fractional diffusion'',
due to its extraordinary capability for describing anomalously slow diffusion processes arising in a wide range of
practical applications in physics, engineering and biology. At a microscopic level, it can be derived from continuous time
random walk with a heavy-tailed waiting time distribution (with a divergent mean) in the sense that the probability
density function of the walker appearing at time $t>0$ and spatial location $x\in\mathbb{R}^d$ satisfies a differential
equation of the form \eqref{eqn:fde}
(in the whole space $\mathbb{R}^d$). The model \eqref{eqn:fde} has been successfully employed in describing
many practical applications, e.g., diffusion of charge carriers in amorphous photoconductors, diffusion in
fractal domains \cite{Nigmatulin:1986}, ion transport in column experiments \cite{HatanoHatano:1998}, and subsurface
flow \cite{AdamsGelhar:1992}. We refer interested readers to the comprehensive reviews \cite{MetzlerKlafter:2000,MetzlerJeon:2014} for
physical motivations of the mathematical model and long lists of successful applications.

The ISPs of interest are to determine some information of the source $F$ from the measurement on
a subboundary $\omega\times\{\ell\}\subset\partial\Omega $ of the domain $\Omega$. Note that the boundary measurement is
insufficient to uniquely determine a general source $F$ (see, e.g., \cite[Section 1.3.1]{KSXY}), and additional
assumptions have to be imposed on the source $F$ in order to restore unique recovery. Often it is formulated as recovering
either spatial or temporal component of the source $F(x,t)$. In this work, the source $F$ is assumed to be of the form
\begin{equation}\label{eqn:cond-f}
  F(x,t) = f(x',t)R(x,t).
\end{equation}
The condition \eqref{eqn:cond-f} can be interpreted as that an unknown source $f(x',t)$ depends only on the
depth variable $x'$ and $t$ in the case of $d=2$, which corresponds to a layer structure, and on the planar location
$(x_1,x_2)$ and $t$ but not on the depth in the case of $d=3$, which can be a good approximation if the domain $\Omega$
is very thin in the direction of $x_3$. Note that it arises also naturally in linearizing the inverse potential problem,
where the potential coefficient $q$ depends only $x'$ and $t$ \cite{Isakov:1990}. We investigate the following two
inverse problems: (i) \textbf{ISPn} is to recover $f(x',t)$ from the boundary observation $u|_{\omega\times\{\ell\}\times(0,T)}$
for $m=1$ in \eqref{eqn:fde} and (ii) \textbf{ISPd} is to recover $f(x',t)$ from the flux measurement $\partial_{x_d}u|_{\omega\times
\{\ell\}\times(0,T)}$ for $m=0$ in \eqref{eqn:fde} (i.e., n and d refer to the Neumann and Dirichlet boundary condition,
respectively, on the subboundary $\omega\times\{\ell\}\times(0,T)$ in the direct problem \eqref{eqn:fde}).

This work is devoted to the theoretical analysis and numerical reconstruction of ISPn and ISPd. In Theorem \ref{thm:uniqueness}, we prove a well-posedness
result for ISPn in $L^2(0,T;L^2(\omega))$. This is achieved by combining the technique developed in \cite{KianYamamoto:2019ip}, improved regularity
estimates and with a novel perturbation argument from \cite{JinLiZhou:2019}. Further, in Theorem \ref{thm:stability}, we
establish a conditional stability result under additional regularity condition on $f(x',t)$ for ISPd. To the best of our knowledge, this is
the first work rigorously analyzing ISPs of recovering a space-time dependent source component in a subdiffusion model with
time-dependent coefficients. The main technical challenges in the study include the nonlocality of the time-fractional
derivative $\partial_t^\alpha u$ and the time-dependence of the operator $\mathcal{A}(t)$. The nonlocality essentially
limits the solution regularity pickup (see, e.g., \cite{SakamotoYamamoto:2011} and \cite[Chapter 6]{Jin:2021}), and thus sharp
regularity estimates for incompatible problem data are needed, which
is especially delicate due to limited smoothness of the domain $\Omega$. This is achieved in Proposition \ref{prop:reg-improved}
by using a refined regularity pickup from \cite[Lemma 2.4]{GaitanKian:2013}, exploiting the cylindrical structure of
the domain $\Omega$. The time dependence of the elliptic operator $\mathcal{A}(t)$ precludes the application of the standard
separation of variable technique that has been predominant in existing studies. This challenge is overcome by a perturbation
argument and maximal $L^p$ regularity for time-fractional problems, which plays an important role in the analysis of ISPd.
In Section \ref{sec:numer}, we derive the adjoint problem for computing the gradient of a quadratic misfit
functional and analyze the regularity of the adjoint variable. Further, we describe the conjugate gradient algorithm for
recovering $f(x',t)$, and provide extensive numerical experiments to illustrate the feasibility of the recovery. The
well-posedness, conditional stability and reconstruction algorithm represent the main contributions of this work.

Last we situate this work in the existing literature. ISPs of recovering a part of information of the source $F$ in a
subdiffusion model from the lateral or terminal data represent an important class of applied inverse problems, and have
been extensively studied in the past decade. Most of the works devoted to this problem have been stated for sources
$F(x, t) = p(t)q(x)$ and can be divided into three groups: (i) inverse $t$-source problem of recovering $p(t)$ \cite{SakamotoYamamoto:2011,
WeiLiLi:2016,FujishiroKian:2016,LiZhang:2020}, (ii) inverse $x$-source problem of recovering $q(x)$ \cite{SakamotoYamamoto:2011b,
ZhangXu:2011,JiangLiLiuYamamoto:2017,RundellZhang:2018}, and (iii) simultaneous inversion of spatial and temporal
components \cite{KSXY,RundellZhang:2020,LiZhang:2020,KJ}. Within group (i), for example, using the decay property of the
Mittag-Leffler function $E_{\alpha,\beta}(z)$, a two-sided stability result of recovering $p(t)$ was shown in \cite{SakamotoYamamoto:2011}, if the observation $u(x_0,t)$
satisfies $x_0\in \text{supp}(q)$. Within group (ii), the unique recovery of the spatial component $q(x)$ by interior observation was proved
in \cite{JiangLiLiuYamamoto:2017} using Duhamel's principle and unique continuation principle, which also gave an iterative
reconstruction algorithm. All these works in groups (i) and (ii) are concerned with recovering only either $p(t)$ or $q(x)$.
The works in group (iii) are closely to the current work. The work \cite{KSXY} showed the simultaneous recovery of $p$ and
$q$ under suitable assumptions. For a two-dimensional heat equation, Rundell and Zhang \cite{RundellZhang:2020}
proved the unique recovery of both $p$ and $q$ in a semi-discrete setting (i.e., the temporal component $p(t)$ is piecewise
constant) from sparse observation on the boundary $\partial\Omega\times (0,T)$. Li and Zhang \cite{LiZhang:2020} extended
the analysis to the time-fractional model in two-dimension, and established the uniqueness of recovering the unknown spatial component $q(x)$,
the time mesh and the fractional order $\alpha$ simultaneously from sparse data on the boundary $\partial\Omega\times(0,T)$. We refer interested readers to the reviews
\cite{JinRundell:2015,LiuLiYamamoto:2019} for further pointers to theoretical and numerical results. See also the work of
\cite{KJ} for the unique recovery of a general source $F$ from the full knowledge of the solution of problem \eqref{eqn:fde},
with $\mathcal{A}$ independent of $t$, on $\Omega\times(T_1, T)$, with $T_1\in (0, T)$.
Kian and Yamamoto \cite{KianYamamoto:2019ip} proved a first uniqueness and stability results for the ISPs of recovering $f(x',t)$ of the
subdiffusion model in a cylindrical domain. (See also Isakov \cite{Isakov:1990} for relevant results for the standard
parabolic problem in the half space.) The analysis in \cite{KianYamamoto:2019ip} relies on some representation of solutions
by mean of $E_{\alpha,\beta}(z)$ which are unavailable for elliptic operators with time-dependent coefficients. This work
extends the results in \cite{KianYamamoto:2019ip} to the case of the time-dependent diffusion coefficients, and further,
by exploiting the maximal $L^p$ regularity, we substantially relax the regularity requirement on $f(x',t)$ for
conditional stability.

Inverse problems for subdiffusion with time-dependent coefficients have been scarcely studied so far, due to a lack of
mathematical tools, when compared with the time-independent counterpart. The only work we are aware of on an ISP
with a time-dependent elliptic operator is \cite{Slodicka:2020}, which showed the unique recovery of a spatial component
from terminal measurement using an energy argument, which seems nontrivial to extend to the case $f(x',t)$. See also
the works \cite{Zhang:2016} for recovering a time-dependent factor in the diffusion coefficient $a(t)$, where the special
structure does allow applying the establish separation of variable technique. Thus, the theoretical analysis for ISPs in the case
of time-dependent coefficients remains challenging. This work presents one promising approach to overcome the challenge
(i.e., perturbation argument), and in particular it allows establishing the stable recovery.

The rest of the paper is organized as follows. In Section \ref{sec:prelim}, we state the
assumptions and preliminary estimates. Then in
Sections \ref{sec:stability} and \ref{sec:stability2}, we prove the well-posedness of ISPn and
conditional stability of ISPd, respectively. In Section \ref{sec:numer},
we describe a numerical algorithm for recovering $f(x',t)$ for both ISPs, and provide several
numerical experiments to showcase the feasibility of the recovery. Throughout, the notation $c$ denotes
a generic constant which may change at each occurrence, but it is always independent of the unknown
source $f(x',t)$ or the associated solution $u$. For a bivariate function $g(x,t)$ or $g(x',t)$, we
often abbreviate it to $g(t)$ as a vector-valued function by suppressing the dependence on the spatial
variable.

\section{Preliminaries: assumptions and basic estimates}\label{sec:prelim}

Now we collect several preliminary results. For $m=0,1$, we define two realizations $A(t)$ and $\tilde A(t)$
in $L^2(\Omega)$ of the elliptic operator $\mathcal{A}$, with their domains respectively given by
\begin{align*}
   D(A(t))&=\{v\in H_0^1(\Omega): \mathcal{A}(t)v\in L^2(\Omega)\},\\ 
   D(\tilde A(t)) &= \{v\in H^1(\Omega): v|_{\partial\omega\times(-\ell,\ell)}=0,\mathcal{A}(t)v\in L^2(\Omega),
   \partial_{x_d}^mv|_{x_d= \ell}=0, \partial_{x_d}v|_{x_d=- \ell}=0\},
\end{align*}
and let $A_*=A(t_*)$ and $\tilde A_*=\tilde A(t_*)$ for any $t_*\in[0,T]$. Note that we abuse the notation $\tilde A(t)$ for both $m=0$ and
$m=1$, which will be clear from the context. For any $s\geq0$, $A_*^s$ and $\tilde A_*^s$ denote the
fractional power of $A_*$ and $\tilde A_*$ via spectral decomposition, and the associated graph
norms by $\|\cdot\|_{D(A_*^s)}$ and $\|\cdot\|_{D(\tilde A_*^s)}$, respectively.
Let $E_*(t)$ and $\tilde E_*(t)$ be the solution operators
corresponding to the source $F$, associated with the elliptic operators $A_*$ and $\tilde A_*$,
respectively, defined by \cite[Section 3.1]{JinLiZhou:2020}
\begin{align}
E_*(t):=\frac{1}{2\pi {\rm i}}\int_{\Gamma_{\theta,\delta }}e^{zt} (z^\alpha+A_*)^{-1}\, \d z  ~~\text{and}~~
&\tilde E_*(t):=\frac{1}{2\pi {\rm i}}\int_{\Gamma_{\theta,\delta}}e^{zt}  (z^\alpha+ \tilde A_*)^{-1}\, \d z , \label{eqn:op}
\end{align}
with the contour $\Gamma_{\theta,\delta}\subset \mathbb{C}$ (oriented with an increasing imaginary part) given by
\begin{equation*}
  \Gamma_{\theta,\delta}=\left\{z\in \mathbb{C}: |z|=\delta, |\arg z|\le \theta\right\}\cup
\{z\in \mathbb{C}: z=\rho e^{\pm\mathrm{i}\theta}, \rho\ge \delta\} .
\end{equation*}
Throughout, we fix $\theta \in(\frac{\pi}{2},\pi)$ so that $z^{\alpha} \in \Sigma_{\alpha\theta}$
for $z\in\Sigma_{\theta}:=\{z\in\mathbb{C}\backslash\{0\}: |{\rm arg}(z)|\leq\theta\}.$ Further, we employ the
operator $\tilde S_*(t)$ (corresponding to the initial data) defined by
\begin{equation*}
  \tilde S_*(t):=\frac{1}{2\pi {\rm i}}\int_{\Gamma_{\theta,\delta}}e^{zt}  z^{\alpha-1}(z^\alpha+ \tilde A_*)^{-1}\, \d z.
\end{equation*}
Then it is known that \cite[(3.8)]{JinLiZhou:2020}
\begin{equation}\label{eqn:diff-S-E}
  \frac{\d}{\d t}\tilde S_*(t) = - \tilde A_*\tilde E_*(t).
\end{equation}

The next lemma summarizes the smoothing properties of $E_*(t)$, $\tilde E_*(t)$ and
$\tilde S_*(t)$. The notation $\|\cdot\|$ denotes the operator norm on $L^2(\Omega)$.
\begin{lemma}[{\cite[Lemma 1]{JinLiZhou:2020}}]\label{lem:smoothing}
For any $\beta\in[0,1]$, there hold for any $t\in (0,T]$
\begin{align*}
   t^{1+\alpha(\beta-1)}\|A_*^\beta  E_*(t)\| & \leq c \quad\mbox{and}\quad
   t^{1+\alpha(\beta-1)}\|\tilde A_*^\beta \tilde E_*(t)\| + t^{1+\alpha}\|\tilde A_*^2\tilde E_*(t)\| + t^{\beta\alpha}\|\tilde A_*^\beta \tilde S_*(t)\| \leq c.
\end{align*}
\end{lemma}

Throughout, we make the following assumption on the diffusion coefficient matrix $a$.
The regularity $a\in C^1([0,T];C^1(\overline{\Omega};\mathbb{R}^{d\times d}))
\cap C([0,T];C^3(\overline{\Omega};\mathbb{R}^{d\times d}))$ is sufficient
for Lemma \ref{lem:perturb}. (ii) is a structural condition to enable unique recovery.
The notation $\cdot$ and $|\cdot|$ denote standard Euclidean inner
product and norm, respectively, on $\mathbb{R}^d$.
\begin{assumption}\label{ass:a}
The coefficient $q\in C^1([0, T];L^\infty(\Omega)) \cap L^\infty(0, T;W^{2,\infty}(\Omega))$, and
the symmetric diffusion coefficient matrix $a\in C^1([0,T];C^1(\overline{\Omega};\mathbb{R}^{d\times d}))
\cap C([0,T];C^3(\overline{\Omega};\mathbb{R}^{d\times d}))$ satisfies the following conditions.
\begin{itemize}
\item[{\rm(i)}]  There exists $\lambda\in(0,1)$ such that for any $(x,t)\in \overline{\Omega}\times[0,T]$,
\begin{equation*}
    \lambda|\xi|^2\leq  a(x,t)\xi\cdot\xi \leq \lambda^{-1}|\xi|^2,\quad \forall \xi\in \mathbb{R}^d.
\end{equation*}
\item[{\rm(ii)}] $a_{jd}(x',\pm\ell,t)=0$, $x'\in\omega$ and $j=1,\ldots,d-1$, and $\partial_{x_d}a_{ij}(t)=0$, for $i,j=1,\ldots,d-1$.
\end{itemize}
\end{assumption}

Note that the cylindrical domain $\Omega=\omega\times(-\ell,\ell)$ is only Lipschitz continuous.
Thus, some extra assumptions on the domain and the coefficient matrix $a$ are needed in order
to guarantee high-order Sobolev regularity of the elliptic operator $\mathcal{A}(t)$ with
suitable boundary conditions. In the analysis, we need the following elliptic regularity
pickup: (i) and (ii) are sufficient for the analysis in Sections \ref{sec:stability} and \ref{sec:stability2},
respectively. (i) holds under the assumption that the domain $\omega$ is convex and
\begin{equation}\label{eqn:conda}
  a_{id}=0,\,\,\partial_{x_j}a_{dd}=0,\,\, \partial_{x_d}a_{ij}=0,\quad i,j\in\{1,\ldots,d-1\}.
\end{equation}
Indeed, if $\omega$ is convex, then $\Omega$
is convex and the desired assertion follows from \cite[Theorems 3.2.1.2 and 3.2.1.3]{Grisvard:1985}.
This can be verified using a separation of variable argument \cite[Lemma 2.4]{GaitanKian:2013}.
Besides the condition \eqref{eqn:conda}, if the domain $\omega$ is of class ${C}^4$,
the separation of variable argument similar to \cite[Lemma 2.4]{GaitanKian:2013} implies
Assumption $\rm \widetilde H$ in Definition \ref{ass:reg}(ii).
\begin{definition}\label{ass:reg}
\begin{itemize}
\item[{\rm(i)}]
A tuple $(\Omega,\mathcal{A}(t))$ is said to satisfy Assumption {\rm H}$mn$, $m,n=0,1$: if for any $t\in[0,T]$ and any $f\in L^2(\Omega)$, the
following boundary value problem
\begin{equation*}
  \left\{\begin{aligned}
    \mathcal{A}(t) v &= f,\quad\mbox{in }\Omega,\\
    v & = 0, \quad \mbox{on } \partial\omega\times(-\ell,\ell),\\
    \partial_{x_d}^mv(x',\ell)& = 0, \quad \mbox{on }\omega,\\
    \partial_{x_d}^nv(x',-\ell)& = 0, \quad \mbox{on }\omega,
  \end{aligned}\right.
\end{equation*}
admits a unique solution $v\in H^2(\Omega)$ such that
$$\|v\|_{H^2(\Omega)}\leq c(\mathcal{A},m,n,\Omega)\|f\|_{L^2(\Omega)}.$$
\item[{\rm(ii)}] A tuple $(\Omega,\mathcal{A}(t))$ is said to satisfy Assumption $\widetilde{\rm H}$, if for all $v\in H^{\max(1,s)}(\Omega)$
satisfying $\mathcal{A}(t)v\in H^s(\Omega)$, $s\in[0,2]$, there holds $v\in H^{2+s}(\Omega)$ and
$$\|v\|_{H^{2+s}(\Omega)}\leq c(\mathcal{A},s,\Omega)(\|\mathcal{A}v\|_{H^s(\Omega)}+\|v\|_{H^s(\Omega)}).$$
\end{itemize}
\end{definition}

The following perturbation estimates are useful.
\begin{lemma}\label{lem:perturb}
Under Assumptions \ref{ass:a}(i) and {\rm H00 / H01 / H11}, for any $t,s\in[0,T]$ and  {$\beta\in[0,1]$}, there hold
\begin{align}
  \|A_*^\beta(I-A(t)^{-1}A(s)) v\|_{L^2(\Omega)}&\leq c|t-s|\|A_*^\beta v\|_{L^2(\Omega)},
  \quad \forall v\in D(A_*^{\beta}),\label{eqn:perturb-1}\\
   \|\tilde A_*^\beta(I-\tilde A(t)^{-1}\tilde A(s))v\|_{L^2(\Omega)}&\leq c|t-s|\|\tilde A_*^\beta v\|_{L^2(\Omega)},
   \quad \forall v\in D(\tilde A_*^{\beta}).\label{eqn:perturb-2}
\end{align}
\end{lemma}
\begin{proof}
For the operator $A(t)$, the case $\beta=0$ is contained in  \cite[Corollary 3.1]{JinLiZhou:2019}. To
show the estimate for $\beta=1$, fix $t,s\in[0,T]$, $v\in D(A_*)$. From Assumption {\rm H}$00$, we deduce
$D(A_*)=H^1_0(\Omega)\cap H^2(\Omega)=D(A(t))=D(A(s))$, i.e., $v\in D(A(t))$ and $v\in D(A(s))$. Moreover, applying again
Assumption {\rm H}$00$, we get
\begin{align*}
  &\|A_*(I-A(t)^{-1}A(s)) v\|_{L^2(\Omega)}\leq c\|(I-A(t)^{-1}A(s)) v\|_{H^2(\Omega)}\\
 \leq &c\|A(t)(I-A(t)^{-1}A(s)) v\|_{L^2(\Omega)}=c\|A(t) v-A(s)v\|_{L^2(\Omega)},
\end{align*}
with $c>0$ a constant independent of $t$ and $s$. Combining this estimate with \cite[eq. (2.6)]{JinLiZhou:2019} and
Assumption {\rm H}$00$, we obtain \eqref{eqn:perturb-1} for $\beta=1$ and $q\equiv0$. We can extend this result to $q\not\equiv0$, since for
$q\in C^1([0,T];L^2(\Omega))$, the mean value theorem implies
\begin{equation*}
  \|q(t)v-q(s)v\|_{L^2(\Omega)}\leq \|\partial_tq\|_{L^\infty(0,T;L^\infty(\Omega))}|t-s|\|v\|_{L^2(\Omega)}\leq c|t-s|\|A_*^\beta v\|_{L^2(\Omega)}.
\end{equation*}
The case $\beta\in(0,1)$ follows by interpolation. The proof of the estimate \eqref{eqn:perturb-2} is identical under
Assumption H01 / H11.
\end{proof}

Below we need Bochner-Sobolev spaces $W^{s,p}(0,T;X)$, for a UMD space $X$ (see \cite{Hytonen:2016}
for the definition of UMD spaces, which include Sobolev spaces $W^{s,p}(\Omega)$ with $1<p<\infty$).
For any $s\ge 0$ and $1\le p< \infty$, we denote by $W^{s,p}(0,T;X)$ the space of functions $v:(0,T)\rightarrow X$,
with the norm defined by complex interpolation. Equivalently, the space is equipped with the quotient norm
\begin{align*}
\|v\|_{W^{s,p}(0,T;X)}&:= \inf_{\widetilde v}\|\widetilde v\|_{W^{s,p}({\mathbb R};X)}
:= \inf_{\widetilde v} \|  \mathcal{F}^{-1}[ (1+|\xi|^2)^{\frac{s}{2}} \mathcal{F}[\widetilde v](\xi) ]\|_{L^p(\mathbb{R};X)}
\end{align*}
where the infimum is taken over all possible $\widetilde v$ that extend $v$ from $(0,T)$ to ${\mathbb R}$,
and $\mathcal{F}$ denotes the Fourier transform (and $\mathcal{F}^{-1}$ being its inverse).
The following norm equivalence result will be used extensively.

\begin{lemma}\label{lem:norm-equiv}
Let $\alpha\in(0,1)$ and $p\in [1,\infty)$ with $\alpha p>1$. If $v(0)=0$ and $ \partial_t^\alpha v \in L^p(0,T;X)$, then  $v \in {W^{\alpha,p}(0,T;X)}$ and
\begin{align*}
 \|v\|_{W^{\alpha,p}(0,T;X)} \le c \|  \partial_t^\alpha v \|_{L^p(0,T;X)}.
 \end{align*}
Meanwhile, if $v(0)=0$, $ v \in  W^{\alpha,p}(0,T;X)$, then $\partial_t^\alpha v \in L^p(0,T; X)$
and
\begin{align*}
\|  \partial_t^\alpha v \|_{L^p(0,T;X)} \le c  \|v\|_{W^{\alpha,p}(0,T;X)}.
\end{align*}
\end{lemma}

\begin{proof}
Let $g = \partial_t^\alpha v \in L^p(0,T;X)$. Then
$v(t) = \frac{1}{\Gamma(\alpha)} \int_0^t (t-s)^{\alpha-1} g(s)\, \d s.$
This and Young's convolution inequality imply (cf., e.g., \cite[Theorem 2.2]{Jin:2021})
$$\|  v \|_{L^p(0,T;X)} \le c \| g \|_{L^p(0,T;X)}.$$
Let $\tilde g$ be the extension of $g$ from $L^p(0,T;X)$ to $L^p(\mathbb{R};X)$
by zero, i.e., $\tilde g(t) = 0$ for $t\in(-\infty,0)\cup(T,\infty)$ and
$\tilde g(t) = g(t)$ for $t\in(0,T)$.
Then let  $$\widetilde v (t) = \frac{1}{\Gamma(\alpha)} \int_{-\infty}^t (t-s)^{\alpha-1} \widetilde g(s)\, \d s,$$
which satisfies
$$\widetilde g(t) = \frac{1}{\Gamma(1-\alpha)}\frac{\d}{\d t}\int_{-\infty}^t (t-s)^{-\alpha} \widetilde v(s)\,\d s.$$
Then there holds  \cite[p. 90]{KilbasSrivastavaTrujillo:2006}
$$({\rm i}\xi)^\alpha \mathcal{F} [\widetilde v](\xi) =  \mathcal{F}[\widetilde  g ](\xi),$$
and $\widetilde v$ is an extension of $v$.
Consequently, we have
\begin{align*}
 \|\widetilde v\|_{W^{\alpha,p}(\mathbb{R};X)}
& =  \|  \mathcal{F}^{-1}[ (1+|\xi|^2)^{\frac{\alpha}{2}} \mathcal{F}[\widetilde v](\xi) ]\|_{L^p(\mathbb{R};X)}
 = \|  \mathcal{F}^{-1}[ K(\xi)  (1+(\mathrm{i}\xi)^\alpha)\mathcal{F}[\widetilde v](\xi) ]\|_{L^p(\mathbb{R};X)}
 \end{align*}
 with $K(\xi) = (1+|\xi|^2)^{\frac{\alpha}{2}}(1+ (\mathrm{i}\xi)^\alpha)^{-1}$.
 Note that
 $$\lim_{|\xi|  \rightarrow0^+} |K(\xi)| = 1 \quad\mbox{and}\quad \lim_{|\xi|  \rightarrow \infty} |K(\xi)| =1,$$
 so $K(\xi)$ is uniformly bounded. Similarly,
\begin{align*}
\xi\frac{\d}{\d\xi} K(\xi)
=& \frac{\alpha |\xi|^2}{1+|\xi|^2} (1+|\xi|^2)^{\frac{\alpha}{2}}(1+(\mathrm{i}\xi)^\alpha)^{-1} +\alpha  (1+|\xi|^2)^{\frac{\alpha}{2}}(1+(\mathrm{i}\xi)^\alpha)^{-2} (\mathrm{i}\xi)^\alpha
\end{align*}
is also bounded. Therefore, vector-valued Mikhlin multiplier theorem (see, e.g. \cite{Bourgain:1986} or \cite[Proposition 3]{Zimmermann:1989}) indicates
that $K(\xi)$ is a Fourier multiplier, and hence
\begin{align*}
&\|v\|_{W^{\alpha,p}(0,T;X)} \le  \|\widetilde v\|_{W^{\alpha,p}(\mathbb{R};X)}
 \le c \|  \mathcal{F}^{-1}[  (1+(\mathrm{i}\xi)^\alpha)\mathcal{F}[\widetilde v](\xi) ]\|_{L^p(\mathbb{R};X)}\\
\le& c \|  \widetilde v \|_{L^p(\mathbb{R};X)} + c \| g  \|_{L^p(\mathbb{R};X)}\le  c \| g  \|_{L^p(\mathbb{R};X)}
= c \| g  \|_{L^p(0,T;X)} = c \| \partial_t^\alpha u \|_{L^p(0,T;X)}.
 \end{align*}
To prove the second assertion, let $v\in C^\infty([0,T];X)$ with $v(0)=0$, and
we extend $v$ from $(0,T)$ to a function $\widetilde v\in W^{\alpha,p}(\mathbb{R};X)$ satisfying
$\widetilde v(t) = 0$ for all $t\leq0$ and
\begin{align} \label{eqn:fbound-02}
\|\widetilde v\|_{W^{\alpha,p}(\mathbb{R};X)} \le c \|v\|_{W^{\alpha,p}(0,T;X)} .
 \end{align}
Then it is direct that
$$ {_{-\infty}\partial_t^\alpha} v(t):=  \frac{1}{\Gamma(\alpha)}\frac{\d}{\d t}\int_{-\infty}^t (t-s)^{-\alpha} \widetilde v(s)\,\d s = \partial_t^\alpha v(t),\quad \forall t\in(0,T),$$
and
\begin{align*}
 &\quad \| \partial_t^\alpha v \|_{L^p(0,T;X)}  = \| {_{-\infty}\partial_t^\alpha} \widetilde v \|_{L^p(0,T;X)} \le  \| {_{-\infty}\partial_t^\alpha} \widetilde v \|_{L^p(\mathbb{R};X)}\\
& = \|  \mathcal{F}^{-1} (\mathrm{i}\xi)^\alpha \mathcal{F}[\widetilde v](\xi) \|_{L^p(\mathbb{R};X)}
 = \|  \mathcal{F}^{-1} K_2(\xi)(1+|\xi|^2)^{\frac{\alpha}{2}} \mathcal{F}[\widetilde v](\xi) \|_{L^p(\mathbb{R};X)}
 \end{align*}
 with $K_2(\xi) = |\xi|^\alpha (1+|\xi|^2)^{-\frac{\alpha}{2}} $. Note that both
 $|K_2(\xi)|$ and $|\xi\frac{\d}{\d\xi}K_2(\xi)| $ are uniformly bounded,
 hence it is a  Fourier multiplier. Then we have
 \begin{align*}
 \| \partial_t^\alpha v \|_{L^p(0,T;X)}&\le \|\mathcal{F}^{-1}(1+|\xi|^2)^{\frac{\alpha}{2}} \mathcal{F}[\widetilde v](\xi) \|_{L^p(\mathbb{R};X)}
 \le c\|\widetilde v\|_{W^{\alpha,p}(\mathbb{R};X)}.
 \end{align*}
 This together with \eqref{eqn:fbound-02} and the density of $C^\infty([0,T];X)$ in ${W^{\alpha,p}(0,T;X)}$
 leads to the second assertion.
\end{proof}

We need the following Gronwall's inequality (see, e.g., \cite{YeGaoDing:2007}, \cite[Exercise 3, p. 190]{Henry:1981}
or \cite[Theorem 4.2]{Jin:2021}).
\begin{lemma}\label{lem:Gronwall}
Let $c,r>0$ and $y,a\in L^1(0,T)$ be nonnegative functions satisfying
\begin{equation*}
  y(t) \leq a(t) + c\int_0^t(t-s)^{r-1}y(s)\d s,\quad t\in(0,T).
\end{equation*}
Then there exists $c=c(r,T)>0$ such that
\begin{align*}
  y(t) \leq a(t) + c\int_0^t (t-s)^{r-1} a(s)\d s,\quad t\in (0,T).
\end{align*}
\end{lemma}

\section{Well-posedness for ISPn}\label{sec:stability}

This section is devoted ISPn, i.e. recovering the source component $f(x',t)$ in problem \eqref{eqn:fde} with $m=1$
from $u|_{\omega\times\{\ell\}\times(0,T)}$. The direct problem is given by
\begin{equation}\label{eqn:fde1}
  \left\{\begin{aligned}
    \partial_t^\alpha u + \mathcal{A}(t) u &= F,\quad \mbox{in }\Omega\times (0,T),\\
    u(x,0) & = 0,\quad \mbox{in }\Omega,\\
    \partial_{x_d}u(x',\pm\ell,t) & = 0, \quad \mbox{on }\omega\times(0,T),\\
   u(x,t) & = 0,\quad \mbox{on } \partial\omega\times(-\ell,\ell) \times (0,T).
  \end{aligned}  \right.
\end{equation}
Subdiffusion with time-dependent coefficients has recently been studied in
\cite{KubicaYamamoto:2018,JinLiZhou:2019,JinLiZhou:2020}, where well-posedness and
several regularity estimates have been established. Our description largely follows
the approach developed in \cite{JinLiZhou:2019,JinLiZhou:2020}.
Throughout, for the prefactor $R(x,t)$ in the source $F$, we make the following assumption.
\begin{assumption}\label{ass:R}
The function $R\in L^\infty(\Omega\times(0,T))$ satisfies $\partial_{x_d}R\in L^\infty (\Omega\times(0,T))$ and that there exists $c_R>0$ such that
$|R(x',\ell,t)|\geq c_R$ for any $(x',t)\in \omega\times(0,T)$.
\end{assumption}

Now we give several regularity estimates for the direct problem \eqref{eqn:fde1}. First we derive a representation of the solution $u$.
The key step is to reformulate problem \eqref{eqn:fde1} into
\begin{equation*}
  \partial_t^\alpha u(t) + \tilde A_* u(t) = F(t)+ (\tilde A_* - \tilde A(t))u(t),\quad \forall  t\in(0,T].
\end{equation*}
According to \cite{JinLiZhou:2019,JinLiZhou:2020}, problem
\eqref{eqn:fde1}  has a unique solution $u$ which satisfies
\begin{equation*}
  u(t) = \int_0^t \tilde E_*(t-s)F(s)\d s + \int_0^t\tilde E_*(t-s) (\tilde A_*-\tilde A(s))u(s)\d s.
\end{equation*}
By setting $t$ to $t_*$, we can use Lemma \ref{lem:perturb}
to estimate the second integral, which involves the crucial perturbation term.

The next result collects \textit{a priori} estimates on the solution $u$ to problem \eqref{eqn:fde1}.
\begin{lemma}\label{lem:reg0}
Let Assumption \ref{ass:reg}(i) hold. Then the solution $u$ to problem \eqref{eqn:fde1} satisfies
\begin{align}\label{eqn:uf}
  \|u(t)\|_{H^1(\Omega)} & \leq c\int_0^t (t-s)^{\frac\alpha2-1}\|F(s)\|_{L^2(\Omega)}\d s,\quad \forall t\in (0,T],
\end{align}
and also the following maximal $L^p$ regularity
\begin{align}\label{eqn:uf-2}
\|\partial_t^\alpha u  \|_{L^p(0,T;L^2(\Omega))}  + \|  u  \|_{L^p(0,T;D(\tilde A_*))}  \le c \|  F  \|_{L^p(0,T;L^2(\Omega))},\quad \forall 1<p<\infty.
\end{align}
\end{lemma}
\begin{proof}
The estimate \eqref{eqn:uf} can be found in \cite[Theorem 2]{JinLiZhou:2020} (with $k=0$),
and \eqref{eqn:uf-2} in \cite[Theorem 2.1]{JinLiZhou:2019}.
\end{proof}

Further, we denote by $u_f$ the solution of problem \eqref{eqn:fde1} to explicitly indicate its
dependence on $f$. First we show that the inverse problem is indeed ill-posed on the space $L^2(0,T;L^2(\omega))$.
\begin{corollary}\label{cor:compact}
Under Assumptions \ref{ass:a}(i) and \ref{ass:R}, the map $f\mapsto u_f|_{L^2(0,T;L^2(\omega))}$ is linear and compact on $L^2(0,T;L^2(\omega))$.
\end{corollary}
\begin{proof}
The linearity is obvious. The compactness is direct from Lemma \ref{lem:reg0}.
In fact, by the maximal $L^p$ regularity in Lemmas \ref{lem:reg0} and \ref{lem:norm-equiv} and Assumption \ref{ass:R}, we have
\begin{equation*}
   \|  u_f \|_{W^{\alpha,2}(0,T;L^2(\Omega))}  +   \|  u_f \|_{L^2(0,T;H^2(\Omega))} \le c \|  f \|_{L^2(0,T;L^2(\omega))}.
\end{equation*}
Thus, $u_f\in W^{\alpha,2}(0,T;L^2(\Omega))\cap L^2(0,T;H^2(\Omega))$. Meanwhile, by
interpolation, the space $W^{\alpha,2}(0,T;L^2(\Omega))\cap L^2(0,T;H^2(\Omega))$
embeds compactly into $L^2(0,T;H^1(\Omega))$ \cite[Theorem 5.2]{Amann:2000},
which, by the trace theorem, embeds continuously into $L^2(0,T;L^2(\omega))$.
Thus the map $f\mapsto u_f|_{\omega\times\{\ell\}\times(0,T)}$ is compact on $L^2(0,T;L^2(\omega))$.
\end{proof}

Let $w = \partial_{x_d} u_f$. Then $w$ satisfies
\begin{equation}\label{eqn:fde-w}
  \left\{\begin{aligned}
    \partial_t^\alpha w + \mathcal{A}(t) w &= -\partial_{x_d}\mathcal{A}(t)u_f(t) + \partial_{x_d}F(t),\quad \mbox{in }\Omega\times (0,T),\\
    w(0) & = 0,\quad \mbox{in }\Omega,\\
    w & = 0, \quad \mbox{on }\partial\Omega \times(0,T).
  \end{aligned}  \right.
\end{equation}
By applying the perturbation argument and using the operator $A(t)$,
the solution $w$ to problem \eqref{eqn:fde-w} can be represented by
\begin{equation}\label{eqn:w-rep}
  w(t) = \int_0^t E_*(t-s)\big(-\partial_{x_d}A(s)u_f(s)+\partial_{x_d}F(s)\big)\d s + \int_0^tE_*(t-s)(A_*-A(s))w(s)\d s.
\end{equation}
Noting the definition $w=\partial_{x_d}u_f$ and the condition $\partial_{x_d}a_{ij}(t)=0$ for $j,j=1,\ldots,d-1$ from
Assumption \ref{ass:a}(ii), we deduce
\begin{align*}
  -(\partial_{x_d} \mathcal{A}(t))u &= \partial_{x_d}(\partial_{x_d}a_{dd}(t)\partial_{x_d}u_f) + \sum_{j=1}^{d-1}\big[\partial_{x_j}(\partial_{x_d}a_{jd}(t)\partial_{x_d}u_f) +
  \partial_{x_d}(\partial_{x_d}a_{jd}(t)\partial_{x_j}u_f)\big]\\
    &\quad - \partial_{x_d}q(t)u_f := B_1(t) w + B_2(t) u_f,
\end{align*}
where the (time-dependent) operators $B_1(t)$ and $B_2(t)$ are respectively given by
\begin{align*}
  B_1(t) w &: = \partial_{x_d}a_{dd}(t)\partial_{x_d} w + 2\sum_{j=1}^{d-1}\partial_{x_d}a_{jd}(t)\partial_{x_j}w + \sum_{j=1}^{d-1}\partial_{x_j}\partial_{x_d}a_{jd}(t)w,\\
  B_2(t) u &: = \sum_{j=1}^d (\partial_{x_d}^2a_{jd}(t))\partial_{x_j}u - \partial_{x_d}q(t)u.
\end{align*}
{Note that Assumption \ref{ass:a}(ii) allows eliminating the cross terms $\partial_{x_i}(\partial_{x_d}
a_{ij}(t)\partial_{x_j}u)$, $i,j=1,\ldots,d-1$, which plays a central role in the analysis below, and without
this, the argument does not work.}

The next result gives useful bounds on $w:=\partial_{x_d}u_f$.
\begin{lemma}\label{lem:reg-w}
Let Assumptions \ref{ass:a}, $\rm H00$, $\rm H11$ and \ref{ass:R} be fulfilled. Then there exists a
unique solution $w\in L^2(0,T;H^2(\Omega))$ with $\mathcal{A}(t)w(t),\partial_t^\alpha w\in L^2(0,T;L^2(\Omega))$
to problem \eqref{eqn:fde-w}, and for any $\beta\in [1,2)$,
\begin{equation*}
  \|w(t)\|_{H^\beta(\Omega)}\leq c\int_0^t(t-s)^{(1-\frac\beta2)\alpha-1}\|f(s)\|_{L^2(\omega)}\d s,\quad t\in (0,T),
\end{equation*}
where the constant $c$ depends only on $R$, $\mathcal{A}$, $\beta$ and $T$.
\end{lemma}
\begin{proof}
By the estimate \eqref{eqn:uf-2} with $p=2$, we have $-\partial_{x_d}\mathcal{A}(t)u_f, ~ \partial_{x_d}F \in L^2(0,T;L^2(\Omega))$.
Then Lemma \ref{lem:reg0} shows that problem \eqref{eqn:fde-w} has a unique solution
$w\in L^2(0,T;H^2(\Omega))$ with $\mathcal{A}(t)w(t),\partial_t^\alpha w\in L^2(0,T;L^2(\Omega))$. Next, we prove the $H^\beta(\Omega)$
bound on $w(t)$. We define the operators $K_1: L^2(0,T;H^1(\Omega))\to L^2(0,T;H^1(\Omega))$ and $K_2: L^2(0,T;L^2(\omega))\to L^2(0,T;H^1(\Omega))$,
respectively, by
\begin{align*}
  K_1v (t) &= \int_0^t E_*(t-s)B_1v(s)\d s,\\
  K_2f (t) & = \int_0^tE_*(t-s)B_2u_f(s)\d s + \int_0^tE_*(t-s)\partial_{x_d}R(s)f(s)\d s.
\end{align*}
By Lemma \ref{lem:smoothing}, we have
\begin{align}
  \|K_1v(t_*)\|_{H^{\beta}(\Omega)} &\leq c\int_0^{t_*}(t_*-s)^{(1-\frac{\beta}{2})\alpha-1}\|B_1v(s)\|_{L^2(\Omega)}\d s\nonumber\\
   &\leq c\int_0^{t_*}(t_*-s)^{(1-\frac\beta2)\alpha-1}\|v(s)\|_{H^1(\Omega)}\d s.\label{eqn:K1}
\end{align}
Similarly, by Lemma \ref{lem:smoothing}, under Assumption \ref{ass:R}, we have
\begin{align}
  \|K_2f(t_*)\|_{H^{\beta}(\Omega)} &\leq c\int_0^{t_*}({t_*}-s)^{(1-\frac\beta2)\alpha-1}\|f(s)\|_{L^2(\omega)}\d s\nonumber\\
   &\quad + c\int_0^{t_*} ({t_*}-s)^{(1-\frac{\beta}{2})\alpha-1}\|u_f(s)\|_{H^1(\Omega)}\d s.\label{eqn:K2-0}
\end{align}
Meanwhile, under Assumption \ref{ass:R} and the estimate \eqref{eqn:uf}, we deduce
\begin{equation*}
  \|u_f(t)\|_{H^1(\Omega)}\leq c\int_0^t(t-s)^{\frac{\alpha}{2}-1}\|f(s)\|_{L^2(\omega)}\d s.
\end{equation*}
Consequently,
\begin{align*}
  \int_0^{t_*} ({t_*}-s)^{(1-\frac{\beta}{2})\alpha-1}\|u_f(s)\|_{H^1(\Omega)}\d s\leq &c\int_0^{t_*} ({t_*}-s)^{(1-\frac{\beta}{2})\alpha-1}\int_0^s(s-\xi)^{\frac\alpha2-1}\|f(\xi)\|_{L^2(\omega)}\d \xi\d s \\
   =& c\int_0^{t_*} \|f(\xi)\|_{L^2(\omega)}\int_\xi^{t_*} (t_*-s)^{(1-\frac\beta2)\alpha-1}(s-\xi)^{\frac\alpha2-1}\d s\d \xi\\
  \leq& cT^\frac\alpha2 \int_0^{t_*}(t_*-s)^{(1-\frac\beta2)\alpha-1}\|f(s)\|_{L^2(\omega)}\d s.
\end{align*}
This and the estimate \eqref{eqn:K2-0} imply
\begin{equation}\label{eqn:K2}
  \|K_2f(t_*)\|_{H^{\beta}(\Omega)} \leq c_T\int_0^{t_*}({t_*}-s)^{(1-\frac\beta2)\alpha-1}\|f(s)\|_{L^2(\omega)}\d s.
\end{equation}
Next, by Lemmas \ref{lem:smoothing} and \ref{lem:perturb}, \eqref{eqn:perturb-1} and Assumption $\rm H00$, we have
\begin{align*}
  &\quad \Big\|\int_0^{t_*}E_*({t_*}-s)(A_*-A(s))w(s)\d s\Big\|_{H^{\beta}(\Omega)} \\
  &\leq c\int_0^{t_*} \|A_*E_*(t_*-s)\|\|A_*^\frac\beta2(I-A_*^{-1}A(s))w(s)\|_{L^2(\Omega)}\d s\\
     & \leq c\int_0^{t_*}(t_*-s)^{-1}(t_*-s) \|w(s)\|_{H^{\beta}(\Omega)}\d s = c\int_0^{t_*}\|w(s)\|_{H^{\beta}(\Omega)}\d s.
\end{align*}
This estimate, \eqref{eqn:K1}, \eqref{eqn:K2} and the solution representation \eqref{eqn:w-rep} lead to
\begin{equation*}
  \|w(t_*)\|_{H^\beta(\Omega)} \leq c\int_0^{t_*}\|w(s)\|_{H^\beta(\Omega)}\d s + c\int_0^{t_*}(t_*-s)^{(1-\frac{\beta}{2})\alpha-1}\|f(s)\|_{L^2(\omega)}\d s.
\end{equation*}
This and Gronwall's inequality in Lemma \ref{lem:Gronwall}
imply the desired $H^\beta(\Omega)$ bound. This completes the proof.
\end{proof}

Now we can state a well-posedness result for ISPn. Note that below we use the notation $L^2(0,T;L^2(\omega))$ and $L^2(\omega\times(0,T))$
interchangeably since they are isomorphic by Fubini--Tonelli theorem.
\begin{theorem}\label{thm:uniqueness}
Let Assumptions \ref{ass:a}, $\rm H00$, $\rm H11$ and \ref{ass:R} be fulfilled. Then for any $f\in L^2(\omega\times(0,T))$, the solution
$u$ of problem \eqref{eqn:fde1} satisfies $u\in H^3(-\ell,\ell; L^2(\omega\times(0,T)))$, $\partial_t^\alpha u$,
$\mathcal{A}(t)u \in H^1(-\ell,\ell;L^2(\omega\times(0,T)))$. Thus, the map
\begin{equation}\label{eqn:h}
  h: (x',t)\mapsto \frac{[\partial_t^\alpha u+(\mathcal{A}(t)+a_{dd}(t)\partial_{x_d}^2 )u](x',\ell,t)}{R(x',\ell,t)}\in L^2(\omega\times(0,T))
\end{equation}
is well defined, and further, there exists a bounded linear operator $\mathcal{H}:L^2(0,T;L^2(\omega))\to L^2(0,T;L^2(\omega))$ such that $f$ solves
\begin{equation}\label{eqn:Hf}
  h=f+\mathcal{H}f,
\end{equation}
which is well-posed on $L^2(\omega\times(0,T))$. Finally, for every pair $(h,f)\in L^2(\omega\times(0,T))\times L^2(\omega\times(0,T))$ satisfying \eqref{eqn:Hf},
the solution $u$ of problem \eqref{eqn:fde1} satisfies \eqref{eqn:h}.
\end{theorem}
\begin{proof}
By Lemmas \ref{lem:reg0} and \ref{lem:reg-w}, problem
\eqref{eqn:fde} has a solution $u_f\in L^2(0,T;H^1(\Omega))$, with
$\mathcal{A}(t)u_f,\partial_t^\alpha u_f\in L^2(0,T;L^2(\Omega))$ and $w=\partial_{x_d}u_f\in L^2(0,T;H^2(\Omega))$,
$\partial_t^\alpha \partial_{x_d}u_f,\ \mathcal{A}(t)\partial_{x_d}u_f\in L^2(0,T;L^2(\Omega))$. Hence,
\begin{align*}
  x_d&\mapsto u_f(\cdot,x_d,\cdot) \in H^3(-\ell,\ell;L^2(\omega\times(0,T))) \cap H^1(-\ell,\ell;L^2(0,T;H^2(\omega))),\\
  x_d&\mapsto \partial_t^\alpha u_f(\cdot,x_d,\cdot) \in H^1(-\ell,\ell;L^2(\omega\times(0,T))).
\end{align*}
By the trace theorem, {we can restrict $\partial_{x_d}w=\partial_{x_d}^2u_f$,
$\partial_t^\alpha u_f$ and $\mathcal{A}(t)u_f$, to the boundary $x_d=\ell$.}
Thus the governing equation in problem \eqref{eqn:fde1} implies that for $(x',t)\in \omega\times(0,T)$,
\begin{equation}\label{eqn:deriv-H}
\begin{split}
  &\quad a_{dd}(t)\partial_{x_d}w(x',\ell,t)  = a_{dd}(t)\partial_{x_d}^2u_f(x',\ell,t)\\
  &= [\partial_t^\alpha u_f +(\mathcal{A}(t)+a_{dd}(t)\partial_{x_d}^2)u_f](x',\ell,t)-R(x',\ell,t)f(x',t)\\
   & = R(x',\ell,t)[h(x',t)-f(x',t)],
\end{split}
\end{equation}
with the function $h(x',t)$ given by \eqref{eqn:h}. Let the operator $\mathcal{H}:L^2(\omega\times(0,T))\to L^2(\omega\times(0,T))$ be defined by
\begin{align*}
[\mathcal{H} \phi] (x',t)= \frac{a_{dd}(t)\partial_{x_d}^2 u_\phi(x',\ell,t)}{R(x',\ell,t)},
\end{align*}
where $u_\phi(x',x_d,t)$ denotes the solution to problem \eqref{eqn:fde1} with $F= \phi R$.
Then it follows from \eqref{eqn:deriv-H} that $f$ is the solution to
\begin{equation*}
   {h =  f + \mathcal{H} f.}
\end{equation*}
Moreover, by Lemma \ref{lem:reg-w}, trace inequality and the defining identity $w=\partial_{x_d}u_f$,
we deduce that for any $\beta\in(\frac34,1)$
\begin{align*}
&\quad\|  f(t) \|_{L^2(\omega)} \le \|  h(t) \|_{L^2(\omega)} + \|  \mathcal{H} f \|_{L^2(\omega)} \le  \|  h(t) \|_{L^2(\omega)} + c \| \partial_{x_d}^2u_f(\cdot,\ell,t)  \|_{L^2(\omega)} \\
&\le  \|  h(t) \|_{L^2(\omega)} + c\| w(t) \|_{H^{2\beta}(\Omega)}
\le  \|  h(t) \|_{L^2(\omega)} + c \int_0^t (t-s)^{(1-\beta)\alpha-1} \| f(s) \|_{L^2(\omega)} \,\d s.
\end{align*}
This and the standard Gronwall's inequality in Lemma \ref{lem:Gronwall} yield
$$\|f(t)\|_{L^2(\omega)}\leq \|h(t)\|_{L^2(\omega)} + c\int_0^t(t-s)^{(1-\beta)\alpha-1}\|h(s)\|_{L^2(\omega)}\d s,$$
which together with Young's inequality directly implies
\begin{equation*}
  \|f\|_{L^2(0,T;L^2(\omega))}\leq c\|h\|_{L^2(0,T;L^2(\omega))}.
\end{equation*}
This shows the well-posedness of equation \eqref{eqn:Hf} and the recovery of $f$ from the data $h$.
Last, fix $(h_1,f)\in L^2(0,T;L^2(\omega))\times L^2(0,T;L^2(\omega))$ satisfying \eqref{eqn:Hf}
with $h=h_1$ and consider $u\in L^2(0,T;H^1(\Omega))$ solving problem \eqref{eqn:fde} with $F=hR$. The preceding
argument shows that one can define $h\in L^2(0,T;L^2(\omega))$ given by \eqref{eqn:h} and $f$ solves
\eqref{eqn:Hf}. This implies {$h_1=f+\mathcal{H}f=h$.} Therefore, we have $h=h_1$, and this completes
the proof of the theorem.
\end{proof}

\begin{remark}
Theorem \ref{thm:uniqueness} actually gives a reconstruction algorithm for recovering $f$ for ISPn, if
the given data $g^\dag(x',t)=u_{f^\dag}(x',\ell,t)$ is sufficiently accurate so that the derivatives
$\partial_t^\alpha g^\dag$ and $\mathcal{A}(t)g^\dag$ in \eqref{eqn:h} can be evaluated accurately. For
noisy data $g^\delta$, one can proceed in two steps: first suitably mollify the data $g^\delta$ so that
the mollified data is smooth, and then apply the fixed point iteration.
\end{remark}
\section{Conditional stability for ISPd}\label{sec:stability2}

In this section, we establish a conditional stability result for ISPd, i.e., recovering $f(x',t)$
in problem \eqref{eqn:fde} with $m=0$ from the lateral flux observation $\partial_{x_d}u|_{\omega
\times\{\ell\}\times(0,T)}$. The direct problem is given by
\begin{equation}\label{eqn:fde0}
\left\{\begin{aligned}
    \partial_t^\alpha u +\mathcal{A}(t)u &=  fR, \quad \mbox{in }\Omega\times (0,T),\\
    u & =0, \quad \mbox{on } \partial\omega\times(-\ell,\ell)\times (0,T),\\
    u(\cdot,\ell,\cdot) & =0 , \quad \mbox{on }\omega\times(0,T),\\
    \partial_{x_d}u(\cdot,-\ell,\cdot)& = 0,\quad \mbox{on }\omega\times(0,T),\\
    u(0) & = 0, \quad \mbox{in }\omega\times(-\ell,\ell).
  \end{aligned}\right.
\end{equation}

Note that the estimates in Lemma \ref{lem:reg0} remain valid for problem \eqref{eqn:fde0}.
The next result gives an improved regularity result, under extra regularity and compatibility
assumptions on the source $F$. This result plays a central role in the stability analysis.
\begin{proposition}\label{prop:reg-improved}
Let Assumptions \ref{ass:a}(i), $\rm H01$ and $\rm\widetilde H$ be fulfilled, $\gamma\in(\frac12,1)$,
and $F\in W^{1,\frac{1}{1+\alpha(1-2\gamma)}}(0,T;L^2(\Omega))\cap L^\frac{1}{\alpha(1-\gamma)}(0,T;H^{2\gamma}(\Omega))$
and $F(0)=0$. Then for any $\beta\in (\frac12,\gamma)$, problem
\eqref{eqn:fde0} has a unique weak solution $u\in L^{\frac{1}{\alpha(1-\gamma)}}(0,T;H^{2+2\beta}(\Omega))
\cap W^{\alpha, \frac{1}{\alpha(1-\gamma)}}(0,T;H^{2\beta}(\Omega))$ with
\begin{align*}
   &\|u\|_{L^{\frac{1}{\alpha(1-\gamma)}}(0,T;H^{2+2\beta}(\Omega))}+ \|u\|_{W^{\alpha, \frac{1}{\alpha(1-\gamma)}}(0,T;H^{2\beta}(\Omega))}\\
  \leq &c\Big(\|F\|_{W^{1,\frac{1}{1+\alpha(1-2\gamma)}}(0,T;L^2(\Omega))}+
  \|F\|_{L^\frac{1}{\alpha(1-\gamma)}(0,T;H^{2\gamma}(\Omega))}\Big).
\end{align*}
\end{proposition}
\begin{proof}
By Sobolev embedding, $F\in L^\infty(0,T;L^2(\Omega))$, and the existence and uniqueness of a
weak solution $u\in L^q(0,T;D(\tilde A_*))$ for all $q\in(1,\infty)$ follows directly from Lemma \ref{lem:reg0} with
\begin{equation}
 \|u\|_{L^q(0,T;D(\tilde A_*))}\leq c\|f\|_{W^{1,\frac{1}{1+\alpha(1-2\gamma)}}(0,T;L^2(\omega))}.
\end{equation}
It suffices to show the claimed regularity. Using the operator $\tilde A(t)$ and
the perturbation argument, since $u(0)=0$, the solution $u$ can be represented by
\begin{equation}\label{eqn:sol-repr-01}
  u(t) = \int_0^t\tilde E_*(s)F(t-s)\d s + \int_0^t\tilde E_*(s)(\tilde A_*-\tilde A(t-s))u(t-s)\d s.
\end{equation}
Then applying $\tilde A_*$ to both sides of the identity, and using the governing equation give
\begin{align*}
  &\partial_t^\alpha u(t)  = -\tilde A_*u(t) +  F(t) + (\tilde A_*-\tilde A(t))u(t)\\
   = & -\int_0^t\tilde A_*\tilde E_*(s)F(t-s)\d s - \int_0^t\tilde A_*\tilde E_*(s)(\tilde A_*-\tilde A(t-s))u(t-s)\d s +  F(t) + (\tilde A_*-\tilde A(t))u(t).
\end{align*}
Now by fixing $t$ at $t_*$ in the identity, applying the identity \eqref{eqn:diff-S-E} and integration by parts formula
to the first integral and noting the condition $F(\cdot,0)=0$ and the fact $\tilde S_*(0)=I$ \cite[Lemma 6.3]{Jin:2021}, we obtain
\begin{align}
  \partial_t^\alpha u(t_*) & = -\int_0^{t_*}\tilde A_*\tilde E_*(s)F(t_*-s)\d s - \int_0^{t_*}\tilde A_*\tilde E_*(s)(\tilde A_*-\tilde A(t_*-s))u(t_*-s)\d s +  F(t_*)\nonumber\\
              & = -\int_0^{t_*}\tilde S_*(s)F'(t_*-s)\d s - \int_0^{t_*}\tilde A_*\tilde E_*(s)(\tilde A_*-\tilde A(t_*-s))u(t_*-s)\d s.\label{eqn:sol-repr-02}
\end{align}
Then it follows from Lemmas \ref{lem:smoothing} and \ref{lem:perturb} that
\begin{align*}
  \|\tilde A_*^{\beta}\partial_t^\alpha u(t_*)&\|_{L^2(\Omega)}
  \leq \int_0^{t_*} \|\tilde A_*^{\beta} \tilde S_*(t_*-s)\| \|F'(s)\|_{L^2(\Omega)} \d s\\
  &\quad + \int_0^{t_*}\|\tilde A_* ^{1+\beta} E_*(t_*-s)\|   \|  \tilde A_* (I- \tilde A_*^{-1} \tilde A(s))u(s)\|_{L^2(\Omega)}\d s\\
  & \leq c\int_0^{t_*} {(t_*-s)^{-\beta\alpha}}\|F'(s)\|_{L^2(\Omega)}\d s + c\int_0^{t_*} {(t_*-s)^{-\beta\alpha-1}(t_*-s)\|\tilde A_*u(s)\|_{L^2(\Omega)}}\d s\\
   & \leq c\int_0^{t_*}{(t_*-s)^{-\beta\alpha}}\| F'(s)\|_{L^2(\Omega)}\d s  + c\int_0^{t_*}{(t_*-s)^{-\beta\alpha}\|\tilde A_*u(s) \|_{L^2(\Omega)}}\d s.
\end{align*}
Note that $t^{-\beta\alpha}\in L^p(0,T)$ for any $p\in (1,\frac{1}{\gamma\alpha}]\subset (1,\frac{1}{\beta\alpha})$,
by the choice $\beta<\gamma$. Now choosing
$r=\frac{1}{\alpha(1-\gamma)}$, $p=\frac{1}{\gamma\alpha}$ and $q=\frac{1}{1+(1-2\gamma)\alpha}$
in the following Young's convolution inequality
\begin{equation}\label{eqn:Young}
  \|f\ast g\|_{L^r(0,T)}\leq \|f\|_{L^p(0,T)} \|g\|_{L^q(0,T)},\quad \forall p,q,r\geq 1\mbox{ with }r^{-1}+1 = p^{-1} + q^{-1},
\end{equation}
we deduce
\begin{align*}
  \|\tilde A_*^{\beta}\partial_t^\alpha u(t_*)\|_{L^{\frac{1}{\alpha(1-\gamma)}}(0,T;L^2(\Omega))}
  &\leq  c\big(\| F'(s)\|_{L^q(0,T;L^2(\Omega))} + \|\tilde A_*u(s) \|_{L^q(0,T;L^2(\Omega))} \big).
\end{align*}
Further, it follows from the representation \eqref{eqn:sol-repr-01} and Lemmas \ref{lem:smoothing} and \ref{lem:perturb} that
\begin{align*}
 \|\tilde A_*^\beta u(t_*)\|_{L^2(\Omega)}&\leq \int_0^{t_*}\|\tilde A_*^\beta \tilde E_*(s)\|\| F(t_*-s)\|_{L^2(\Omega)}\d s \\
   &\qquad + \int_0^{t_*}\|\tilde A_*\tilde E_*(s)\|\|\tilde A_*^\beta (I-\tilde A_*^{-1}\tilde A(t_*-s))u(t_*-s)\|_{L^2(\Omega)}\d s\\
 &\leq c\int_0^{t_*}s^{(1-\beta)\alpha-1}\| F(t_*-s)\|_{L^2(\Omega)}\d s + c\int_0^{t_*}s^{-1}s\|u(t_*-s)\|_{L^2(\Omega)}\d s\\
 & \leq c\|F\|_{L^\infty(0,T;L^2(\Omega))} t_*^{(1-\beta)\alpha} + c \int_0^{t_*}\|\tilde A_*^\beta u(s)\|_{L^2(\Omega)}\d s.
\end{align*}
This and Gronwall's inequality directly imply $\lim_{t\to0^+}\|\tilde A_*^{\beta} u(t)\|_{L^2(\Omega)} = 0$.
Hence, from Lemma \ref{lem:norm-equiv} and Assumption ${\rm H01}$, we deduce $u\in W^{\alpha,\frac{1}{\alpha(1-\gamma)}}(0,T;D(\tilde A_*^{\beta}))
\subset W^{\alpha,\frac{1}{\alpha(1-\gamma)}}(0,T; H^{2\beta}(\Omega))$. Thus, we conclude that for any fixed $t\in(0,T]$, the solution $u$ satisfies
\begin{equation*} 
  \left\{\begin{aligned}
   \mathcal{A}(t) u (t) &= F(t)- \partial_t^\alpha u(t) ,\quad \mbox{in }\Omega\\
    u(x',\ell,t) & = 0, \quad \mbox{on }\omega ,\\
    \partial_{x_d} u(x',-\ell,t) & = 0, \quad \mbox{on }\omega ,\\
   u(x,t) & = 0,\quad \mbox{on } \partial\omega\times(-\ell,\ell) .
  \end{aligned}  \right.
\end{equation*}
Note that for any $t\in(0,T]$, there holds
$  \mathcal{A}(\cdot) u (\cdot)  = F(\cdot)- \partial_t^\alpha u(\cdot) \in L^{\frac{1}{\alpha(1-\gamma)}}(0,T;H^{2\beta}(\Omega))$.
Then by Assumption $\rm \widetilde H$, we obtain $u\in L^{\frac{1}{\alpha(1-\gamma)}}(0,T;H^{2(1+\beta)}(\Omega))$. This completes the proof.
\end{proof}

The conditional stability analysis employs the regularity estimates in Proposition \ref{prop:reg-improved}.
Let $u$ be the solution to problem \eqref{eqn:fde0}, and let $v=\partial_{x_d}u$. Then
$v$ satisfies
\begin{equation}\label{eqn:v}
  \left\{\begin{aligned}
    \partial_t^\alpha v +\mathcal{A}(t)v &= H + f\partial_{x_d}R, \quad \mbox{in }\Omega\times (0,T],\\
    v & =0, \quad \mbox{on } \partial\omega\times(-\ell,\ell)\times (0,T),\\
    v(\cdot,\ell,\cdot) & =\partial_{x_d}u(\cdot,\ell,\cdot) , \quad \mbox{on }\omega\times(0,T],\\
    v(\cdot,-\ell,\cdot)& = 0,\quad \mbox{on }\omega\times(0,T),\\
    v(0) & = 0, \quad \mbox{in }\omega\times(-\ell,\ell),
  \end{aligned}\right.
\end{equation}
with the function $H$ given by
\begin{equation*}
  H(x',x_d,t) = -\partial_{x_d}\mathcal{A}(t) u= \sum_{i,j=1}^d\partial_{x_i}(\partial_{x_d}a_{ij}(t)\partial_{x_j}u) - \partial_{x_d}q(t)u.
\end{equation*}
	
Unlike problem \eqref{eqn:fde-w} in Section \ref{sec:stability}, problem \eqref{eqn:v} involves a nonzero
Dirichlet boundary condition, and thus requires a different analysis. We employ an extension
approach to derive the requisite bound on $v$. For $r\geq 1$, the notation $X_{\alpha,r}$ denotes the
the space $W^{\alpha,r}(0,T;L^2(\omega))\cap L^r(0,T;H^\frac{3}{2}(\omega))$ with the norm
\begin{equation*}
  \|v\|_{X_{\alpha,r}} = \|v\|_{W^{\alpha,r}(0,T;L^2(\omega))} + \|v\|_{L^r(0,T;H^\frac{3}{2}(\omega))}.
\end{equation*}

\begin{assumption}\label{ass:stab}
$\omega\in C^4$,  $f\in W^{1,\frac{1}{1+(1-2\gamma)\alpha}}(0,T;L^2(\omega))
\cap L^{\frac{1}{\alpha(1-\gamma)}}(0,T;H^{2\gamma}(\omega))$, for some $\gamma\in (\frac34,1)$, with $f(0)=0$, $\partial_t R\in L^\infty
(\Omega\times(0,T))$ and $R\in C([0,T]; W^{2,\infty}(\Omega))$.
\end{assumption}

\begin{lemma}\label{lem:reg-v}
Let Assumptions \ref{ass:a}, \ref{ass:R}, ${\rm H00}$, ${\rm H01}$ and \ref{ass:stab} be
fulfilled, and let $u$ be the solution of problem
\eqref{eqn:fde0}. Then for any $\beta\in(\frac34,\gamma)$, the solution $v$ to problem \eqref{eqn:v} satisfies
\begin{align*}
  \|\partial_{x_d}v(\cdot,\ell,t)\|_{L^2(\omega)}\leq&  c\|\partial_{x_d}u(\cdot,\ell,\cdot)\|_{X_{\alpha,\frac{1}{\alpha(1-\gamma)}}} +
   c\int_0^t(t-s)^{\alpha(1-\beta)-1}\|f(\cdot,s)\|_{L^2(\omega)}\d s,
\end{align*}
for any $t\in(0,T)$, where the constant $c$ depend on $R$, $\Omega$, $T$, $\alpha$, $\beta$, $\gamma$ and $\mathcal{A}$.
\end{lemma}
\begin{proof}
Let $r=\frac{1}{\alpha(1-\gamma)}$.
{Since $f\in W^{1,\frac{1}{1+\alpha(1-2\gamma)}}(0,T;L^2(\omega))$ and $\partial_t R\in L^{\infty}(\Omega\times(0,T))$,
we deduce $F=fR\in  W^{1,\frac{1}{1+\alpha(1-2\gamma)}}(0,T;L^2(\Omega))$.
The assumptions $f\in L^r(0,T;H^{2\gamma}(\omega))$ and
$R\in C([0,T];W^{2,\infty}(\Omega))$ imply $F=L^\frac{1}{\alpha(1-\gamma)}(0,T;H^{2\gamma}(\Omega))$.
Further, $f(\cdot, 0)=0$ in $\omega$ indicates $F(\cdot, 0)=0$ in $\Omega$.}
Thus, $F=fR$ satisfies the conditions in Proposition \ref{prop:reg-improved}, and since $\omega\in C^4$,
Assumption $\rm \widetilde H$ holds. By Proposition \ref{prop:reg-improved}, $u\in W^{\alpha,r}(0,T;H^{2\beta}(\Omega))\cap
 L^r(0,T;H^{2(1+\beta)}(\Omega))$, for any $\beta\in (\frac34,\gamma)$, and by the trace theorem, there holds
\begin{equation}\label{eqn:reg-u-3}
  (x',t)\mapsto \partial_{x_d}u(x',\ell,t)\in
  W^{\alpha,r}([0,T];H^{2\beta-\frac{3}{2}}(\omega))\cap
  L^r([0,T];H^{\frac{1}{2}+2\beta}(\omega)),
\end{equation}
and $\partial_{x_d}u(\cdot,-\ell,0)=0$. Next we split the solution $v$ to
problem \eqref{eqn:v} into $v=v_1+v_2$,
with the functions $v_1$ and $v_2$, respectively, solving
\begin{equation*}
  \left\{\begin{aligned}
    \partial_t^\alpha v_1 + \mathcal{A}(t)v_1 & = 0,\quad \mbox{in }\Omega\times(0,T),\\
    v_1 &=0,\quad \mbox{on }\partial\omega\times(-\ell,\ell)\times(0,T),\\
    v_1(\cdot,\ell,\cdot) &= \partial_{x_d} u(\cdot,\ell,\cdot) ,\quad \mbox{on }\in \omega\times(0,T),\\
    v_1(\cdot,-\ell,t) &= 0,\quad \mbox{on } \omega\times(0,T),\\
    v_1(0) & = 0,\quad \mbox{in }\Omega,
  \end{aligned} \right.
\end{equation*}
and
\begin{equation*}
  \left\{\begin{aligned}
    \partial_t^\alpha v_2 + \mathcal{A}(t)v_2 & = H+f\partial_{x_d}R,\quad \mbox{in }\Omega\times(0,T),\\
    v_2 &=0,\quad \mbox{on }\partial\Omega\times(0,T),\\
    v_2(0) & = 0,\quad \mbox{in }\Omega.
  \end{aligned} \right.
\end{equation*}
Next we bound $v_1$ and $v_2$.  To bound $v_1$, we first extend $\partial_{x_d}
u(\cdot,\ell,\cdot)$ from $\omega\times (0,T)$ to $\Omega\times (0,T)$.
Indeed, by the regularity estimate \eqref{eqn:reg-u-3} and using the classical lifting theorem
for Sobolev spaces \cite[Chapter 1, Theorem 9.4]{LionsMagenes:1972} and Assumption H$00$,
we deduce that there exists a function $G\in  L^r(0,T;H^2(\Omega))$ satisfying
\begin{equation}\label{cond1}
-\Delta_x G(x,t)=0,\quad (x,t)\in Q
\end{equation}
\begin{align}\label{cond2}
  &G(x',x_d,t) = \left\{\begin{aligned}
    \partial_{x_d}u(x',\ell,t),& \quad  x_d=\ell, (x',t)\in \omega\times(0,T),\\
    0,  &\quad x_d= -\ell, (x',t)\in\omega\times(0,T),\\
    0, &\quad t=0, (x',x_d)\in \omega\times(-\ell,\ell).
  \end{aligned}\right.
\end{align}
Clearly, $G$ satisfies the following estimate
\begin{equation*}
  \|G\|_{L^r(0,T;H^2(\Omega))} \leq c\|\partial_{x_d}u(\cdot,\ell,\cdot)\|_{L^r(0,T;H^\frac{3}{2}(\omega))}.
\end{equation*}
Moreover, by interpreting $G$ as the solution to  \eqref{cond1}-\eqref{cond2} in the transposition
sense (see e.g. \cite[Chapter 2, Theorem 6.3]{LionsMagenes:1972} or  \cite{Berggren:2004}), $G\in
W^{\alpha,r}(0,T;L^2(\Omega))$ satisfies
\begin{equation*}
  \|G\|_{W^{\alpha,r}(0,T;L^2(\Omega))} \leq c\|\partial_{x_d}u(\cdot,\ell,\cdot)\|_{W^{\alpha,r}(0,T;L^2(\omega))}.
\end{equation*}
Consequently, we have
\begin{align}
   \|G\|_{W^{\alpha,r}(0,T;L^2(\Omega))}+\|G\|_{L^r(0,T;H^2(\Omega))} \leq \|\partial_{x_d}u(\cdot,\ell,\cdot)\|_{X_{\alpha,r}}.\label{eqn:ext-G}
\end{align}
Then we can decompose $v_1$ into $v_1=G+w_1,$
with the function $w_1$ solving
\begin{equation*}
  \left\{\begin{aligned}
    \partial_t^\alpha w_1 +\mathcal{A}(t)w_1 & = F_1,\quad \mbox{in }\Omega\times(0,T),\\
    w_1 &=0,\quad \mbox{on }\partial\Omega\times(0,T),\\
    w_1(0) & = 0,\quad \mbox{in }\omega\times(-\ell,\ell),
  \end{aligned} \right.
\end{equation*}
with $F_1=-\partial_t^\alpha G-\mathcal{A}(t)G.$ Since  $\partial_{x_d}u(x', \ell, 0) =0$ for $x'\in\omega$,
the uniqueness of the solution of problem \eqref{cond1}--\eqref{cond2} implies
$G(\cdot,0) = 0$. Then, direct computation with Lemma \ref{lem:norm-equiv} gives
\begin{align}
  \|F_1\|_{L^r(0,T;L^2(\Omega))}& \leq \|\partial_t^\alpha G\|_{L^r(0,T;L^2(\Omega))}+
  \|\mathcal{A}(t)G\|_{L^r(0,T;L^2(\Omega))}\nonumber\\
    &\leq c(\|G\|_{W^{\alpha,r}(0,T;L^2(\Omega))}+\|G\|_{L^r(0,T;H^2(\Omega))}).\label{eqn:est-F1}
\end{align}
Thus using the operator $A(t)$ and the perturbation argument, we have
\begin{equation*}
  w_1(t_*) = \int_0^{t_*}E_*(t_*-s)F_1(s)\d s + \int_0^{t_*}E_*(t_*-s)(A_*-A(s))w_1(s)\d s.
\end{equation*}
By Lemmas \ref{lem:smoothing} and \ref{lem:perturb},
\begin{align*}
  \|A_*^\beta w_1(t_*)&\|_{L^2(\Omega)} \leq \int_0^{t_*}\|A_*^\beta E_*(t_*-s)\|\|F_1(s)\|_{L^2(\Omega)} \d s \\
  &\quad + \int_0^{t_*}\|A_*E_*(t_*-s)\|\|A_*^\beta(I-A_*^{-1}A(s))w_1(s)\|_{L^2(\Omega)}\d s\\
  &\leq c\int_0^{t_*}(t_*-s)^{(1-\beta)\alpha-1}\|F_1(s)\|_{L^2(\Omega)} \d s
  + c\int_0^{t_*}(t_*-s)^{-1}(t_*-s)\|A_*^\beta w_1(s)\|_{L^2(\Omega)}\d s\\
   &\leq c\|F_1\|_{L^r(0,T;L^2(\Omega))} + \int_0^{t_*}\|A_*^\beta w_1(s)\|_{L^2(\Omega)}\d s.
\end{align*}
It follows from this estimate and Gronwall's inequality that $w_1\in L^\infty(0,T;D(A_*^\beta))$ with
\begin{equation*}
  \|w_1\|_{L^\infty(0,T;D(A_*^\beta)))}\leq c_T\|F_1\|_{L^r(0,T;L^2(\Omega))}.
\end{equation*}
Then by the triangle inequality, \eqref{eqn:est-F1} and Assumption H00,
\begin{align*}
  \|v_1\|_{L^\infty(0,T;D(A_*^\beta))}&\leq \|w_1\|_{L^\infty(0,T;D(A_*^\beta))} + \|G\|_{L^\infty(0,T;D(A_*^\beta))}\\
  &\leq c\big(\|F_1\|_{L^r(0,T;L^2(\Omega))}+\|G\|_{L^\infty(0,T;{H^{2\beta}}(\Omega))}\big).
\end{align*}
Meanwhile,
the condition $\beta\in (\frac34,\gamma)$ and \cite[Theorem 5.2]{Amann:2000} imply the following embedding inequality
\begin{equation*}
  \|w\|_{L^\infty(0,T;H^{2\beta}(\Omega))} \leq c(\|w\|_{W^{\alpha,r}(0,T;L^2(\Omega))} + \|w\|_{L^r(0,T;H^2(\Omega))}).
\end{equation*}
The last two estimates together give
\begin{equation*}
  \|v_1\|_{L^\infty(0,T;D(A_*^\beta))}
  \leq c\big(\|G\|_{W^{\alpha,r}(0,T;L^2(\Omega))}+\|G\|_{L^r(0,T;H^2(\Omega))}\big).
\end{equation*}
This and the estimate \eqref{eqn:ext-G} imply
\begin{equation}\label{eqn:bdd-v1-0}
  \|v_1\|_{L^\infty(0,T;D(A_*^\beta))}\leq c\|\partial_{x_d}u(\cdot,\ell,\cdot)\|_{X_{\alpha,r}}.
\end{equation}
Moreover, by Assumption H00 and the trace inequality
$$\|\partial_{x_d}v_1(\cdot,\ell,\cdot)\|_{L^\infty(0,T;L^2(\omega))}\leq c\|v_1\|_{L^\infty(0,T;D(A_*^\beta))},$$
we get
\begin{equation}\label{eqn:bdd-v1}
  \|\partial_{x_d}v_1(\cdot,\ell,\cdot)\|_{L^\infty(0,T;L^2(\omega))}\leq c \|\partial_{x_d}u(\cdot,\ell,\cdot)\|_{X_{\alpha,r}}.
\end{equation}
Next we bound $v_2$. Note that the solution $v_2(t)$ can be represented by
\begin{equation*}
  v_2(t) = \int_0^{t}E_*(t-s)[H(s)+\partial_{x_d} F(s)]\d s + \int_0^{t}E_*(t-s)(A_*-A(s))v_2(s)\d s.
\end{equation*}
Thus, by Lemmas \ref{lem:smoothing} and \ref{lem:perturb} and Assumption \ref{ass:R}, we get
\begin{align}
  \|A_*^{\beta} v_2(t_*)\|_{L^{2}(\Omega)}\leq& c\int_0^{t_*}(t_*-s)^{\alpha(1-\beta)-1}\big[\|f(s)\|_{L^2(\omega)} +\|H(s)\|_{L^2(\Omega)}\big]\d s\nonumber\\
    &+ c\int_0^{t_*} \|A_*^{\beta} v_2(s)\|_{L^{2}(\Omega)}\d s.\label{eqn:est-v2}
\end{align}
In light of Assumption \ref{ass:a}(ii) and the definition $v=\partial_{x_d}u$, we have
\begin{align*}
   &H(t)  = - \partial_{x_d}\mathcal{A}(t)u\\
  =& \partial_{x_d}a_{dd}(t)\partial_{x_d}^2 u + 2\sum_{j=1}^{d-1}\partial_{x_d}a_{jd}(t)\partial_{x_j}\partial_{x_d} u + \sum_{j=1}^d\partial_{x_j}\partial_{x_d}a_{jd}(t)\partial_{x_d} u \\
   &\quad + \sum_{j=1}^{d-1}\partial_{x_d}^2a_{jd}(t)\partial_{x_j} u-\partial_{x_d}q(t)u\\
  =& \partial_{x_d}a_{dd}(t)\partial_{x_d}v + 2\sum_{j=1}^{d-1}\partial_{x_d}a_{jd}(t)\partial_{x_j}v + \Big(\sum_{j=1}^d\partial_{x_j}\partial_{x_d}a_{jd}(t)\Big)v + \sum_{j=1}^{d-1}\partial_{x_d}^2a_{jd}(t)\partial_{x_j} u - \partial_{x_d}q(t)u,
\end{align*}
from which it directly follows that
\begin{equation*}
  \|H(t)\|_{L^2(\Omega)}\leq c\big(\|v(t)\|_{H^1(\Omega)} +\|u(t)\|_{H^1(\Omega)}\big),\quad t\in (0,T].
\end{equation*}
By Lemma \ref{lem:reg0} with $\beta=\frac12$ (which holds also for problem \eqref{eqn:fde0}) and Assumption \ref{ass:R}, we have
\begin{equation*}
  \|u(t)\|_{H^1(\Omega)}\leq c\int_0^t(t-s)^{\frac{\alpha}{2}-1}\|f(s)\|_{L^2(\omega)}\d s.
\end{equation*}
This and \eqref{eqn:bdd-v1-0} lead to
\begin{align*}
  \|H(t)\|_{L^2(\Omega)}\leq & c\int_0^t(t-s)^{\frac\alpha2-1}\|f(s)\|_{L^2(\omega)}\d s+ c(\|\partial_{x_d}u(\cdot,\ell,\cdot)\|_{X_{\alpha,r}}
  +\|v_2(t)\|_{H^1(\Omega)}),
\end{align*}
which together with \eqref{eqn:est-v2} yields
\begin{align*}
  \|A_*^\beta v_2(t_*)\|_{L^{2}(\Omega)}\leq& c\int_0^{t_*}(t_*-s)^{(1-\beta)\alpha-1}\|f(s)\|_{L^2(\omega)}\d s + c\|\partial_{x_d}u(\cdot,\ell,\cdot)\|_{X_{\alpha,r}}\\
    & + c\int_0^{t_*}(t_*-s)^{(1-\beta)\alpha-1} \|A_*^\beta v_2(s)\|_{L^{2}(\Omega)}\d s.
\end{align*}
This estimate and Gronwall's inequality in Lemma \ref{lem:Gronwall} then imply
\begin{align*}
  \|A_*^\beta v_2(t)\|_{L^2(\Omega)}\leq c\int_0^t(t-s)^{(1-\beta)\alpha-1}\|f(\cdot,s)\|_{L^2(\omega)}\d s+ c\|\partial_{x_d}u(\cdot,\ell,\cdot)\|_{X_{\alpha,r}}. \end{align*}
It follows from this estimate, Assumption H00 and the trace inequality that
\begin{align*}
  \|\partial_{x_d}v_2(\cdot,\ell,t)\|_{L^2(\omega)}
\leq c\int_0^t(t-s)^{\alpha(1-\beta)-1}\|f(s)\|_{L^2(\omega)}\d s + c\|\partial_{x_d}u(\cdot,\ell,\cdot)\|_{X_{\alpha,r}}.
\end{align*}
Finally, combining this bound with the estimate \eqref{eqn:bdd-v1} yields the desired assertion.
\end{proof}

Now we can state a conditional stability result for ISPd.
\begin{theorem}\label{thm:stability}
Let Assumptions \ref{ass:a}, \ref{ass:R}, $\rm H00$, $\rm H01$ and \ref{ass:stab} be fulfilled, and $u$ the solution
of problem \eqref{eqn:fde0}. Then there exists a constant $c$
depending on $R$, $\Omega$, $T$, $\alpha$, $\gamma$, $p$ and $\mathcal{A}$ such that
\begin{equation*}
  \|f\|_{L^\infty(0,T;L^2(\omega))}\leq c\Big(\|\partial_{x_d}u(\cdot,\ell,\cdot)\|_{L^\frac{1}{\alpha(1-\gamma)}(0,T;H^\frac{3}{2}(\omega))}+\|\partial_{x_d}u(\cdot,\ell,\cdot)\|_{W^{\alpha
  ,\frac{1}{\alpha(1-\gamma)}}(0,T;L^2(\omega))}\Big).
\end{equation*}
\end{theorem}
\begin{proof}
First projecting the governing equation in \eqref{eqn:fde0} onto the lateral boundary $\omega\times\{\ell\}\times(0,T)$ and then using the fact
that, for all $(x',t)\in \omega \times(0,T)$, we have $u(x',\ell,t)=0$, and thus for any $(x',t)\in\omega\times(0,T)$,
\begin{equation*}
  f(x',t)R(x',\ell,t) =-\Big[a_{dd}(t)\partial_{x_d}^2u + 2\sum_{j=1}^{d-1}a_{jd}(t)\partial_{x_j}\partial_{x_d}u+\sum_{j=1}^d\partial_{x_j}a_{jd}(t)\partial_{x_d}u\Big](x',\ell,t).
\end{equation*}
This, Assumption \ref{ass:R}, and the definition $v=\partial_{x_d}u$ imply that for all $t\in(0,T)$, there holds
\begin{align*}
  \|f(t)\|_{L^2(\omega)} &\leq c_R^{-1}c\big(\|\partial_{x_d}v(\cdot,\ell,t)\|_{L^2(\omega)}+\|\partial_{x_d}u(\cdot,\ell,t)\|_{H^1(\omega)}\big).
\end{align*}
Under the condition $\gamma\in (\frac34,1)$, the choice $r=\frac{1}{\alpha(1-\gamma)}$, \cite[Theorem 5.2]{Amann:2000} implies
\begin{equation*}
  \|\partial_{x_d}u(\cdot,\ell,t)\|_{L^\infty(0,T;H^1(\omega))}\leq c(\|\partial_{x_d}u(\cdot,\ell,t)\|_{W^{\alpha,r}(0,T;L^2(\omega))}+\|\partial_{x_d}u(\cdot,\ell,t)\|_{L^r(0,T;H^\frac32(\omega))}).
\end{equation*}
The last two estimates and Lemma \ref{lem:reg-v} imply
\begin{align*}
  \|f(t)\|_{L^2(\omega)} &\leq c\|\partial_{x_d}u(\cdot,\ell,\cdot)\|_{X_{\alpha,r}}
    + c\int_0^t(t-s)^{\alpha(1-\beta)-1}\|f(s)\|_{L^2(\omega)}\d s.
\end{align*}
Then Gronwall's inequality in Lemma \ref{lem:Gronwall} implies the desired assertion, completing the proof of the theorem.
\end{proof}

\begin{remark}
Theorem \ref{thm:stability} shows the influence of the fractional order $\alpha$ on the stability: the larger is the order $\alpha$,
the stronger temporal regularity  $\|\partial_{x_d}u(\cdot,\ell,\cdot)\|_{W^{\alpha,\frac{1}{\alpha(1-\gamma)}}(0,T;L^2(\omega))}$ on the
data $u|_{\omega\times\{\ell\}\times (0,T)}$ the stability needs. This agrees with the smoothing
property of the solution operator, and shows also the beneficial influence of anomalous diffusion. 
Theorem \ref{thm:stability} improves the corresponding result in \cite[Theorem 1.4]{KianYamamoto:2019ip} {\rm(}with $\delta>\frac12${\rm)}:
\begin{equation*}
  \|f\|_{L^\infty(0,T;L^2(\omega))}\leq c(\|\partial_{x_d}u(\cdot,\ell,\cdot)\|_{L^\infty(0,T;H^\frac32(\omega))}+\|\partial_{x_d}u(\cdot,\ell,\cdot)\|_{W^{1,\infty}(0,T;H^\delta(\omega))}),
\end{equation*}
This improvement is achieved by the maximal $L^p$ regularity and the suitable interpolation inequality in
fractional Sobolev spaces.
\end{remark}

\begin{remark}
In the spirit of \cite[Corollary 1.5]{KianYamamoto:2019ip}, Theorem \ref{thm:stability} allows proving the stable recovery
of a class of the zeroth order coefficient $q$ from the flux data $\partial_{x_d}u|_{\omega\times\{\ell\}\times(0,T)}$.
This analysis requires the existence of a solution to problem \eqref{eqn:fde0} in $W^{1,\infty}(0,T;W^{1,\infty}(\Omega))\cap L^{\infty}
(0,T;W^{2,\infty}(\Omega))$. The latter can be achieved using the argument of Proposition \ref{prop:reg-improved}, and we
leave the details to future investigation.
\end{remark}

\section{Numerical experiments and discussions}\label{sec:numer}
In this section, we present several numerical experiments to illustrate the feasibility
of recovering the space-time dependent  $f$ from lateral boundary observation.

\subsection{Numerical algorithm}

First we describe a numerical algorithm for recovering $f$ for ISPn (and the algorithm for ISPd is
similar). We employ an iterative regularization technique, which approximately minimizes
\begin{equation}\label{eqn:Tikh}
  J(f) :=\tfrac12\|u_f-g^\delta\|_{L^2(0,T;L^2(\omega))}^2,
\end{equation}
where $u_f$ denotes the solution to the direct problem \eqref{eqn:fde1} with $F=fR$. By Corollary \ref{cor:compact},
the map $u_f:L^2(0,T;L^2(\omega))\to L^2(0,T;L^2(\omega))$ is linear and compact, and thus standard regularization theory
\cite{EnglHankeNeubauer:1996,ItoJin:2015} can be applied to justify the reconstruction technique. In particular,
when equipped with an appropriate stopping criterion, the approximate minimizer obtained by gradient
type methods, e.g., gradient descent and conjugate gradient method, will converge to the exact source component $f^\dag$ as the noise level
tends to zero, and further it will converge at a certain rate dependent of the ``regularity'' of $f^\dag$ (in the
sense of source condition or conditional stability estimates), when equipped with
a suitable stopping criterion.

To (approximately) minimize the functional $J(f)$, we employ the conjugate gradient (CG) method
\cite{AlifanovArtyukhin:1995}. When applying the method, the main computational effort is to compute the gradient, which can be done efficiently
using the adjoint technique. Specifically, let $v$ be the solution to the following adjoint problem
\begin{equation}\label{eqn:fde-adj}
  \left\{\begin{aligned}
    _t\kern-.5em^R\kern-.2em\partial_T^\alpha v + \mathcal{A}(t) v &= 0,\quad (x',x_d,t)\in \Omega\times (0,T),\\
    {_tI_T^{1-\alpha}}v(x,T) & = 0,\quad \mbox{in }\Omega,\\
    \partial_{x_d} v(x',\ell,t) & = u_f-g^\delta, \quad \mbox{on }\omega\times(0,T),\\
    \partial_{x_d} v(x',-\ell,t) & = 0, \quad \mbox{on }\omega\times(0,T),\\
   v(x,t) & = 0,\quad \mbox{on } \partial\omega\times(-\ell,\ell) \times (0,T),
  \end{aligned}  \right.
\end{equation}
where the notation $_tI_T^{1-\alpha} v(t)$ and $_t\kern-.5em^R\kern-.2em\partial_T^\alpha v$
denotes the right-sided Riemann-Liouville fractional integral and derivative of $v$,
defined  respectively by  \cite{KilbasSrivastavaTrujillo:2006}
\begin{equation*}
 _tI_T^{1-\alpha} v(t) = \frac{1}{\Gamma(1-\alpha)}\int_t^T(s-t)^{-\alpha}v(s)\d s\quad \mbox{and}\quad
   _t\kern-.5em^R\kern-.2em\partial_T^\alpha v(t) = -\frac{1}{\Gamma(1-\alpha)}\frac{\d}{\d t}\int_t^T (s-t)^{-\alpha}u(s)\d s.
\end{equation*}

Then we have the following representation of the gradient $J'(f)$ of $J(f)$.
\begin{proposition}
The gradient $J'(f)$ of the functional $J(f)$ is given by
\begin{equation}\label{eqn:grad}
  J'(f) = \int_{-\ell}^\ell Rv\,\,\d x_d,
\end{equation}
where $v$ is the solution to the adjoint problem \eqref{eqn:fde-adj}.
\end{proposition}
\begin{proof}
The derivation follows a standard procedure. The directional derivative $J'(f)[h]$ of
the functional $J$ with respect to $f$ in the direction
$h\in L^2(0,T;L^2(\omega))$ is given by
\begin{equation*}
  J'(f)[h] = (u_h,u_f-g^\delta)_{L^2(0,T;L^2(\omega))},
\end{equation*}
where $u_h$ is the solution to problem \eqref{eqn:fde1} with $h$ in place of $f$ (or the source $F=hR$).
Multiplying the equation for $u_h$ with a test function $\phi(x,t)$ and then integration by parts yield
\begin{equation}\label{eqn:fde-lin}
  \int_0^T\int_\Omega (\phi\,\partial_t^\alpha u_h+ a\nabla u_h\cdot\nabla \phi) \d x\d t = \int_0^T\int_\Omega Rh\phi\d x \d t.
\end{equation}
Meanwhile, the weak formulation for the adjoint solution $v$ is given by
\begin{equation}\label{eqn:fde-adj-weak}
  \int_0^T\int_\Omega (\phi\, {_t\kern-.5em^R\kern-.2em\partial_T^\alpha}v + a\nabla v\cdot\nabla \phi)\d x\d t = \int_0^T\int_\omega (u_f-g^\delta)\phi\d x'\d t.
\end{equation}
Then taking $\phi=v$ in \eqref{eqn:fde-lin} and $\phi=u_h$ in \eqref{eqn:fde-adj-weak}, appealing to the identity
(\cite[p. 76, Lemma 2.7]{KilbasSrivastavaTrujillo:2006} or \cite[Lemma 2.6]{Jin:2021})
$$\int_0^T\int_\Omega v\partial_t^\alpha u_h \d x\d t = \int_0^T\int_\Omega u_h{_t\kern-.5em^R\kern-.2em\partial_T^\alpha} v\d x\d t$$
 (in view of the zero initial / terminal conditions)
and subtracting the two identities give
\begin{equation*}
  \int_0^T\int_\Omega Rhv\d x \d t = \int_0^T\int_\omega (u_f-g^\delta)u_h\d x'\d t.
\end{equation*}
This and the definition of the derivative $J'(f)$ show the desired assertion.
\end{proof}

The next result gives the regularity of the adjoint variable $v$.
\begin{theorem}\label{thm:reg-adj}
Let $g^\delta\in L^2(0,T;L^2(\omega))$, and Assumption ${\rm H11}$ be fulfilled. Then there exists a unique solution
$v\in L^2(0,T;H^{\frac{3}{2}-\epsilon}(\Omega))\cap W^{\frac{3\alpha}{4}-\epsilon,2}(0,T;L^2(\Omega))$, for any
$\epsilon>0$, to the adjoint problem \eqref{eqn:fde-adj}.
\end{theorem}
\begin{proof}
For any fixed $t\in[0,T]$, let $N(t)$ be the Neumann map defined by $\phi=N(t)\psi$, with $\phi$ solving
\begin{equation}\label{eqn:neuman-map}
  \left\{\begin{aligned}
 \mathcal{A}(t) \phi &= 0,\quad \mbox{in }\Omega,\\
    \partial_{x_d} \phi(x',-\ell) & = 0, \quad \mbox{on }\omega,\\
      \partial_{x_d} \phi(x',\ell) & = \psi, \quad \mbox{on }\omega,\\
   \phi(x) & = 0,\quad \mbox{on } \partial\omega\times(-\ell,\ell).
  \end{aligned}  \right.
\end{equation}
It is known that $ \|N(t) \psi\|_{H^{\frac32}(\Omega)} \le c\| \psi  \|_{L^2(\omega)}$ with a range
$\mathcal{R}(N(t)) = D(\tilde A^{\frac34-\epsilon})$, for any small $\epsilon>0$ \cite[Proposition 2.12]{Acquistapace:1991}.
Below we first analyze the case of time-independent coefficients, and then the case of time-dependent
coefficients. Let $\psi=u_f-g^\delta\in L^2(0,T;L^2(\omega))$. \vskip5pt
\noindent\textbf{Case i:} $\mathcal{A}(t) \equiv  \mathcal{A}_*$ and $N(t)\equiv N_*$.
Note that the solution $v$ to problem \eqref{eqn:fde-adj} can be represented by using the
operators $\tilde E(t)$ and $N(t)$
\cite[Section 2.2]{KianYamamoto:2020}:
\begin{equation} \label{eqn:sol-adj}
v(t) =  \int_t^{T} \tilde A_* \tilde E_*(s-t)   N_* \psi(s) \,\d s.
\end{equation}
Thus, by Lemma \ref{lem:smoothing} and Young's inequality,
for any $\theta\in[0,\frac34)$ and $\epsilon\in(0,2-2\theta-\frac12)$, there holds
\begin{equation*}
\begin{split}
 \| v \|_{L^2(0,T; D(\tilde A_*^{\theta}))} &
  \le \int_0^T \|  \tilde A_*^{\theta+\frac14+\frac\epsilon2} E_*(t)  \|\,\d t
  \Big( \int_0^T \|  \tilde A_*^{\frac34-\frac\epsilon2} N_*  \psi(t) \|_{L^2(\Omega)}^2\,\d t \Big)^{\frac12}
  \le c_\epsilon \| \psi\|_{L^2(0,T;L^2(\omega))}.
\end{split}
\end{equation*}
It follows from equation \eqref{eqn:fde-adj} that
\begin{equation*}
  {_t\kern-.5em^R\kern-.2em\partial_T^\alpha} v  +  \tilde A_* ( v-N_*\psi )= 0,
\end{equation*}
which implies ${_t\kern-.5em^R\kern-.2em\partial_T^\alpha} v \in L^2(0,T; D(\tilde A_*^{\theta-1}))$.
Then by a similar argument as in the proof of Lemma \ref{lem:norm-equiv},
we deduce $v \in  W^{\alpha,2}(0,T; D(\tilde A_*^{\theta-1}))$ (see also \cite[Theorem 2.1]{JinLiZhou:control}).
Then by interpolation, we derive $v\in W^{\frac{3\alpha}{4}-\epsilon,2}(0,T;L^2(\Omega))$ for any $\epsilon>0$ \cite[Theorem 5.2]{Amann:2000}.
Further, in view of the identity
\begin{equation*}
  {_tI_T^{1-\alpha}} v (t) = \int_t^T \tilde A_* \tilde S_*(s-t) N_* \psi(s) \,\d s,
\end{equation*}
and Young's inequality, we deduce for any $\theta\in(1-\frac1{2\alpha},\frac34)$
\begin{equation*}
\begin{split}
 &\| {_tI_T^{1-\alpha}} v (t)\|_{L^2(\Omega)} \le  \int_t^T \|\tilde A_*^{1-\theta} \tilde S_*(s-t)\|\, \|\tilde A_*^{\theta} N_* \psi (s)\|_{L^2(\Omega)} \,\d s\\
\le &  c \int_t^T(s-t)^{-(1-\theta)\alpha}\, \|\psi(s)\|_{L^2(\omega)} \,\d s
\le c \Big(\int_t^T(s-t)^{-2(1-\theta)\alpha}\, \,\d s\Big)^{\frac12} \, \|  \psi \|_{L^2(0,T;L^2(\omega))} \\
\le &c (T-t)^{\frac12-(1-\theta)\alpha} \|\psi \|_{L^2(0,T;L^2(\omega))} .
\end{split}
\end{equation*}
Therefore the terminal condition ${_tI_T^{1-\alpha}} v (T) = 0$ holds.\\
\medskip
\noindent\textbf{Case ii:} time-dependent elliptic operator $\mathcal{A}(t)$. We rewrite the adjoint problem \eqref{eqn:fde-adj} as
\begin{equation}\label{eqn:fde-adj-theta}
  {_t\kern-.5em^R\kern-.2em\partial_T^\alpha} v(t) + \tilde A_* (v- N\psi )(t) =   \big(\tilde A_*- \tilde A(t)\big)
(v- N\psi )(t),
\end{equation}
with ${_tI_T^{1-\alpha}} v (T) = 0$. Then the solution $v(t)$ can be represented as
\begin{equation} \label{eqn:sol-adj-t}
v(t) =  \int_t^{T} \tilde A_* \tilde E_*(s-t)  N(t)\psi(s) \,\d s
+  \int_t^{T} \tilde A_* \tilde E_*(s-t)\big(\tilde A_*- \tilde A(s)\big)
(v(s)- N(s)\psi(s))\,\d s.
\end{equation}
By Lemmas \ref{lem:smoothing} and \ref{lem:perturb}, we deduce that for any
$\theta\in[0,\frac34)$ and $\epsilon\in(0,2-2\theta-\frac12)$,
\begin{align}
\| v(t_*) \|_{D(\tilde A^\theta)}\le &c \int_{t_*}^T(s-t_*)^{-(\theta+\frac14+\frac\epsilon2)\alpha}\| \psi(s) \|_{L^2(\omega)}\d s\nonumber\\
    &+ c \int_{t_*}^T \|  v(s) \|_{D(\tilde A^\theta)}\d s +  c \int_{t_*}^T \|\psi (s) \|_{L^2(\omega)} \d s. \label{eqn:est-adj-01}
\end{align}
Then squaring both sides of \eqref{eqn:est-adj-01} and integrating over $[t_0,T]$ lead to
\begin{equation*}
  \|v\|_{L^2(t_0,T;D(\tilde A^\theta))}^2 \le c_\epsilon \|\psi \|_{L^2(0,T; L^2(\omega))}^2
  + c \int_{t_0}^T \|  v \|_{L^2(t,T;D(\tilde A^\theta))}^2  \d t.
\end{equation*}
This together with Gronwall's inequality implies that for any $\theta\in[0,\frac34)$
\begin{equation*}
  \|v\|_{L^2(0,T;D(\tilde A^\theta))} \le c  \| \psi \|_{L^2(0,T;L^2(\omega))}.
\end{equation*}
By \eqref{eqn:fde-adj-theta}, ${_t\kern-0.5em^R\kern-.2em\partial_T^\alpha} v \in L^2(0,T; D(\tilde A_*^{\theta-1}))$, and by Lemma \ref{lem:norm-equiv},
$ v \in  W^{\alpha,2}(0,T; D(\tilde A_*^{\theta-1}))$. Then by interpolation, we derive
$u\in W^{\frac{3\alpha}{4}-\epsilon,2}(0,T;L^2(\Omega))$ for any $\epsilon>0$ \cite[Theorem 5.2]{Amann:2000}.
\end{proof}
\begin{remark}
It follows from Theorem \ref{thm:reg-adj} that for data $g^\delta\in L^2(0,T;L^2(\omega))$,
the gradient $J'(f)$ belongs to $L^2(0,T;H^{\frac32-\epsilon}(\omega))\cap W^{\frac{3\alpha}{4}-\epsilon,2}
(0,T;L^2(\omega))$ for any small $\epsilon>0$, if the factor $R$ is smooth, and further,
$_tI_T^{1-\alpha} J'(f)(x',T)=0$ for $x'\in\omega$ and $J'(f)(x',t)=0$ for $(x',t)\in\partial
\omega\times(0,T)$. These conditions will impact the convergence behavior of the conjugate gradient method,
dependent of the regularity of $f^\dag$.
\end{remark}

Now we can describe the conjugate gradient method \cite{AlifanovArtyukhin:1995} for minimizing  $J$.
The complete procedure is listed in Algorithm \ref{alg:cgm}. In the algorithm, Steps 6-7 compute the
conjugate descent direction, and Step 8 computes the optimal step size using the sensitivity problem.
In general, the algorithm converges within tens of iterations; see the numerical experiments below.
At each iteration, the algorithm involves three forward solves (direct problem, adjoint problem and
sensitivity problem), which represent the main computational effort. For the stopping criterion
at Step 11, we employ the discrepancy principle \cite{Morozov:1966,EnglHankeNeubauer:1996,ItoJin:2015}, i.e.,
\begin{equation}\label{eqn:dp}
 k^* = \arg\min \{k\in\mathbb{N}: \|u_{f^k}-g^\delta\|_{L^2(0,T;L^2(\omega))} \leq c\delta\},
\end{equation}
where $c>1$  and $\delta=\|g^\dag-g^\delta\|_{L^2(0,T;L^2(\omega))}$
is the noise level of the data $g^\delta$.

\begin{algorithm}[hbt!]
  \caption{Conjugate gradient method for minimizing the functional $J$ in \eqref{eqn:Tikh}.\label{alg:cgm}}
  \begin{algorithmic}[1]
    \STATE Initialize $f^0$, and set $k=0$.
    \FOR{$k=0,\ldots,K$}
    \STATE Solve for $u^k$ from problem \eqref{eqn:fde1} with $F=f^kR$, and compute the residual $r^k=u^k|_{\omega\times\ell\times (0,T)}-g^\delta$.
    \STATE Solve for $v^k$ from problem \eqref{eqn:fde-adj} with $r^k$.
    \STATE Compute the gradient $J'(f^k)$  by \eqref{eqn:grad}.
    \STATE Compute the conjugate coefficient $\gamma^k$ by
    \begin{equation*}
      \gamma^k = \left\{\begin{aligned}
      0, &\quad k= 0,\\
      \frac{\|J'(f^k)\|_{L^2(0,T;L^2(\omega))}^2}{\|J'(f^{k-1})\|_{L^2(0,T;L^2(\omega))}^2}
      , &\quad k\geq 1.
      \end{aligned}\right.
    \end{equation*}
    \STATE Compute the conjugate direction $d^k$ by $d^k=-J'(f^k)+\gamma^kd^{k-1}$.
    \STATE Solve for $u_{d^k}$ from problem \eqref{eqn:fde1} with $F=d^kR$.
    \STATE Compute the step size $s^k$ by
    \begin{equation*}
      s^k = -\frac{(u_{d^k},r^k)_{L^2(0,T;L^2(\omega))}}{\|u_{d^k}\|^2_{L^2(0,T;L^2(\omega))}}.
    \end{equation*}
    \STATE Update the source component $f^{k+1}=f^k+s^k d^k$.
    \STATE Check the stopping criterion.
    \ENDFOR
  \end{algorithmic}
\end{algorithm}

Algorithm \ref{alg:cgm} can also be applied to ISPd, by viewing the zero Dirichlet
data on $\omega\times\{\ell\}\times(0,T)$ as the measurement,
and then the measurement $\partial_{x_d}u|_{\omega\times\{\ell\}\times(0,T)}$
as the Neumann data on $\omega\times\{\ell\}\times(0,T)$ for problem \eqref{eqn:fde1}.
However, the discrepancy principle \eqref{eqn:dp} cannot be applied directly, due to
a lack of the noise level for the Dirichlet boundary data.

\subsection{Numerical results and discussions}

Now we present several examples to illustrate the feasibility of recovering $f$.
The domain $\Omega$ is taken to be the unit square $\Omega=(-\frac12,\frac12)^2$, with $\omega=(-\frac12,
\frac12)$, $q\equiv0$, and the final time $T=1$. The direct and adjoint problems are all discretized by the standard
continuous piecewise linear Galerkin method in space, and backward Euler convolution quadrature in time;
see \cite{JinLiZhou:2019} for the error analysis for relevant direct problems and the review
\cite{JinLazarovZhou:2019} for various numerical schemes. The domain $\Omega$ is first divided into $M^2$
small squares each of width $1/M$, and then further divided into triangles by connecting the upper right vertex
with the lower left vertex of each small square to obtain a uniform triangulation. For the inversion step, we
take $M=100$ and $N=1000$. The same spatial and temporal mesh is used for approximating $f$. The factor
$R(x,t)$ is fixed at $R\equiv 1$. The exact data $g^\dag$ on the lateral boundary $\omega\times\{\ell\}
\times(0,T)$ is obtained by solving the direct problem \eqref{eqn:fde} with the exact $f^\dag$ on a finer
mesh. The noisy boundary data $g^\delta$ is generated from the exact data $g^\dag$ by
\begin{equation*}
  g^\delta(x',t) = g^\dag(x',t) + \varepsilon\|g^\dag\|_{L^\infty(\omega\times(0,T))}\xi(x',t),\quad \forall (x',t)\in\omega\times(0,T),
\end{equation*}
where $\xi(x',t)$ follows the standard normal distribution, and $\varepsilon$ denotes the
relative noise level. In Algorithm \ref{alg:cgm}, the maximum number of CG iterations is fixed at $50$,
and the constant $c$ in \eqref{eqn:dp} is taken to be $c=1.01$.
Throughout, we measure the accuracy of a reconstruction $\hat f$ by the $L^2$ error $e(\hat f)$ defined by
\begin{equation*}
e(\hat f) = \|\hat f-f^\dag\|_{L^2(0,T;L^2(\omega))}.
\end{equation*}

\subsubsection{Numerical results for ISPn}
First we illustrate the case of time-independent coefficients.
\begin{example}\label{exam:ind}
$a(x_1,x_2) = 1+\sin(\pi x_1)x_2(1-x_2)$ and $f^\dag(x_1,t) =(\frac14-x_1^2)t(T-t)e^t$.
\end{example}

\begin{figure}[hbt!]
  \centering
  \setlength{\tabcolsep}{0pt}
  \begin{tabular}{ccc}
    \includegraphics[width=0.33\textwidth]{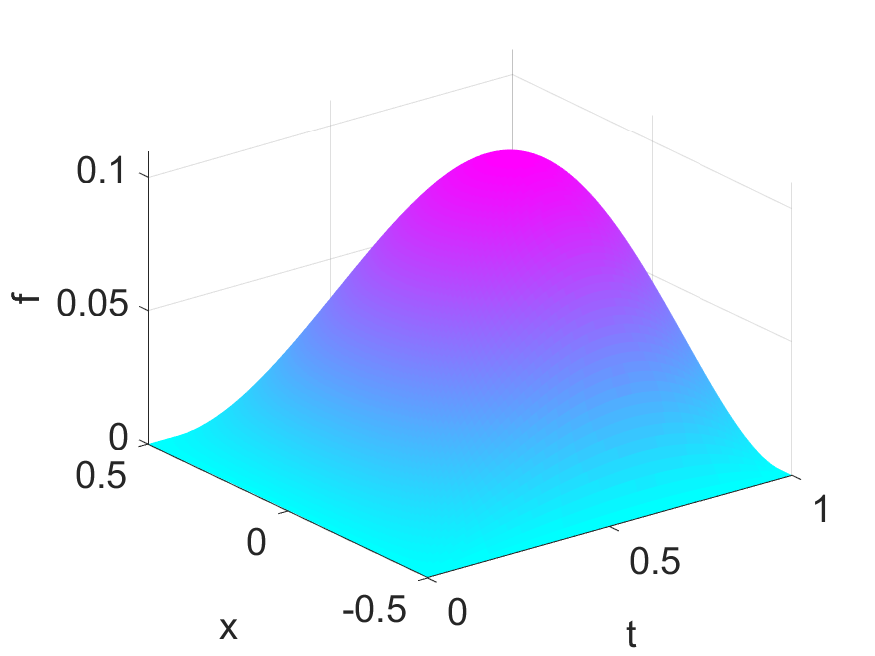} & \includegraphics[width=0.33\textwidth]{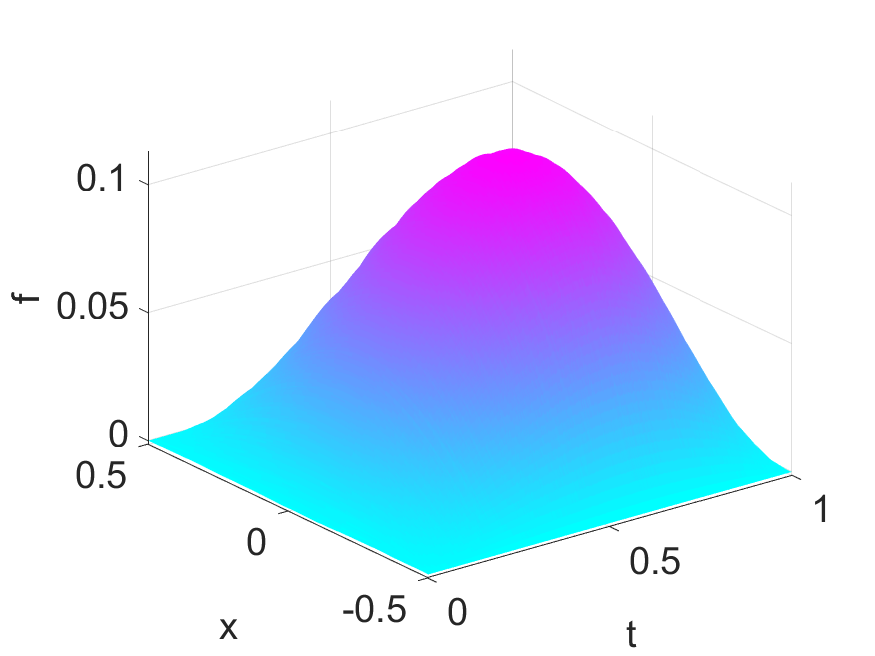} & \includegraphics[width=0.33\textwidth]{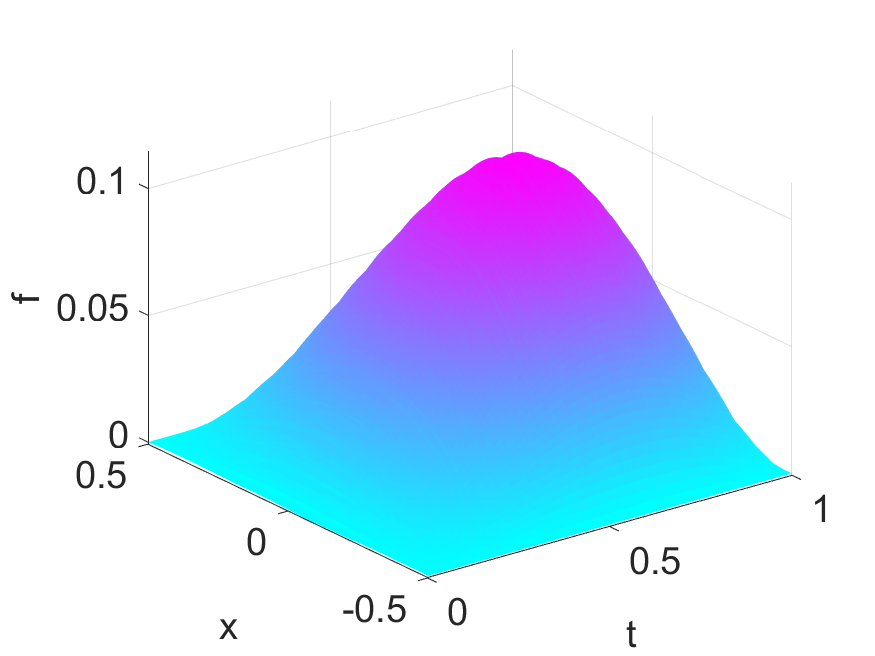}\\
       & \includegraphics[width=0.33\textwidth]{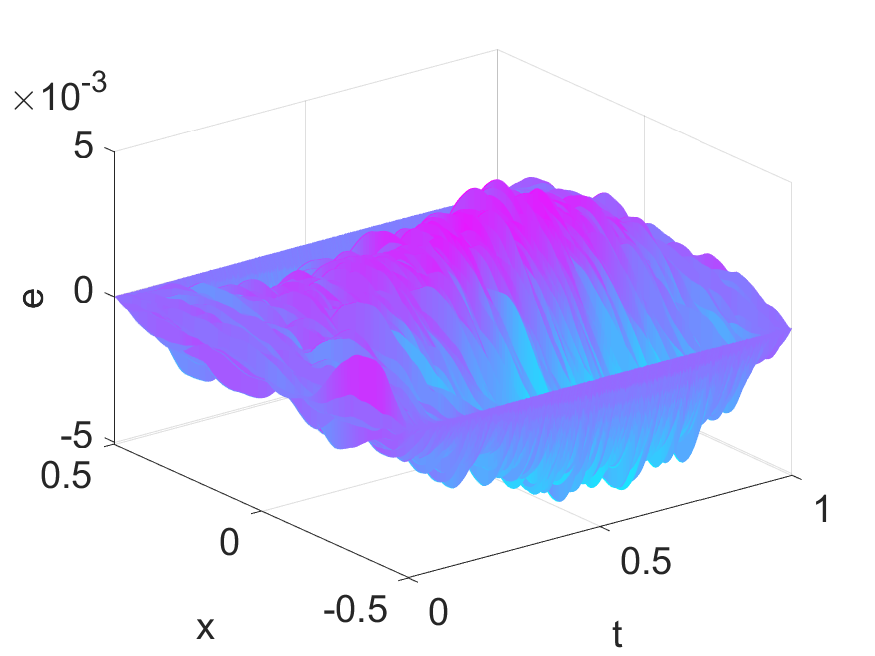} & \includegraphics[width=0.33\textwidth]{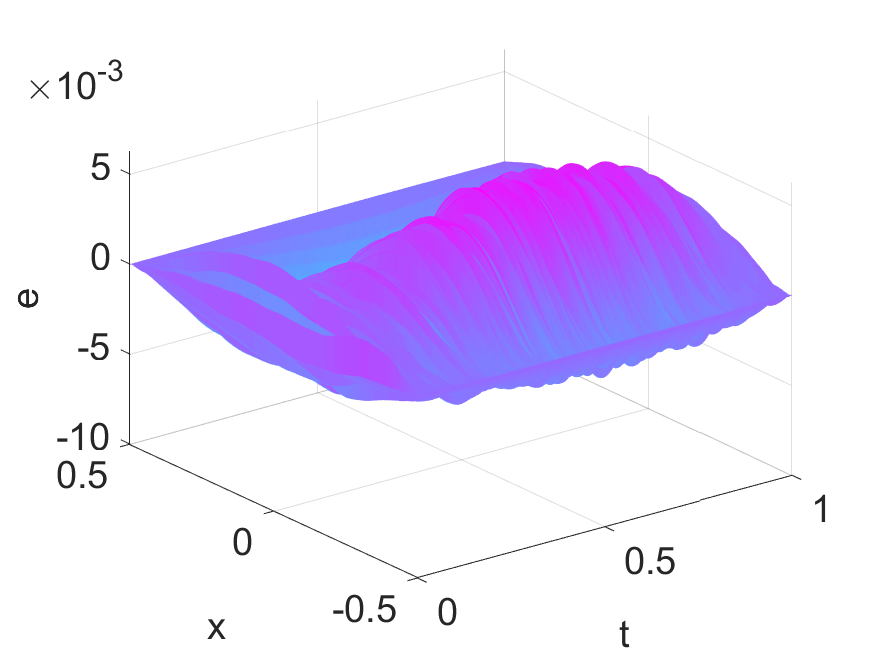}\\
       (a) exact & (b) $\varepsilon$=1e-2 & (c) $\varepsilon$=5e-2
  \end{tabular}
  \caption{Reconstructions and the pointwise errors for Example \ref{exam:ind} with $\varepsilon$=1e-2 and $\varepsilon$=5e-2.\label{fig:ind}}
\end{figure}

The numerical results for Example \ref{exam:ind} are presented in Table \ref{tab:ind}, where the numbers
in the bracket denote the stopping index determined by the discrepancy principle \eqref{eqn:dp}. For
noisy data $g^\delta$, the method reaches convergence within ten iterations, and thus it is fairly efficient.
It is observed that as the relative noise level $\epsilon$ increases from zero to 5e-2, the error $e(\hat f)$ also
increases, whereas the required number of CG iterations decreases. For a fixed noise level $\varepsilon$, the
reconstruction error $e$ tends to decrease with the order $\alpha$, and all the reconstructions are
fairly accurate; see Fig. \ref{fig:ind} for typical reconstructions and the associated pointwise errors
$e=\hat f-f^\dag$ (which slightly abuses the notation $e$). These results clearly show the feasibility of recovering $f$
from the lateral boundary data, corroborating the theoretical results in \cite{KianYamamoto:2019ip}.

\begin{table}[hbt!]
\centering
\caption{The reconstruction errors $e$ for Example \ref{exam:ind}.\label{tab:ind}}
\begin{tabular}{cccccc}
\toprule
$\alpha\backslash\varepsilon$ & 0 & 1e-3 & 5e-3 & 1e-2 & 5e-2\\
\midrule
0.25 & 8.61e-5 (50) & 3.87e-4 (13) & 7.36e-4 (10) & 1.26e-3 (7) & 2.33e-3 (4)\\
0.50 & 4.27e-5 (50) & 3.91e-4 (10) & 6.84e-4 ( 8) & 1.29e-3 (6) & 2.19e-3 (3)\\
0.75 & 8.61e-5 (50) & 3.62e-4 (16) & 5.93e-4 (11) & 8.84e-4 (9) & 1.71e-3 (4)\\
\bottomrule
\end{tabular}
\end{table}

Now we give two examples with time-dependent coefficients. The notation $\chi_S$ denotes the characteristic
function of a set $S$.
\begin{example}\label{exam:dep}
The diffusion coefficient $a$ is given by $a(x_1,x_2,t)=(1+\sin(\pi x_1)x_2(1-x_2))(1+\sin t)$, and consider two different source components.
\begin{itemize}
  \item[{\rm(i)}] $f^\dag (x_1,t) = \sin(x_1+\frac12)\pi t(T-t)e^t$.
  \item[{\rm(ii)}] $f^\dag (x_1,t) = \sin(x_1+\frac12)\pi t(T-t)e^t\chi_{[0,0.7]}(t)$.
\end{itemize}
\end{example}

\begin{figure}[hbt!]
  \centering
  \setlength{\tabcolsep}{0pt}
  \begin{tabular}{ccc}
    \includegraphics[width=0.33\textwidth]{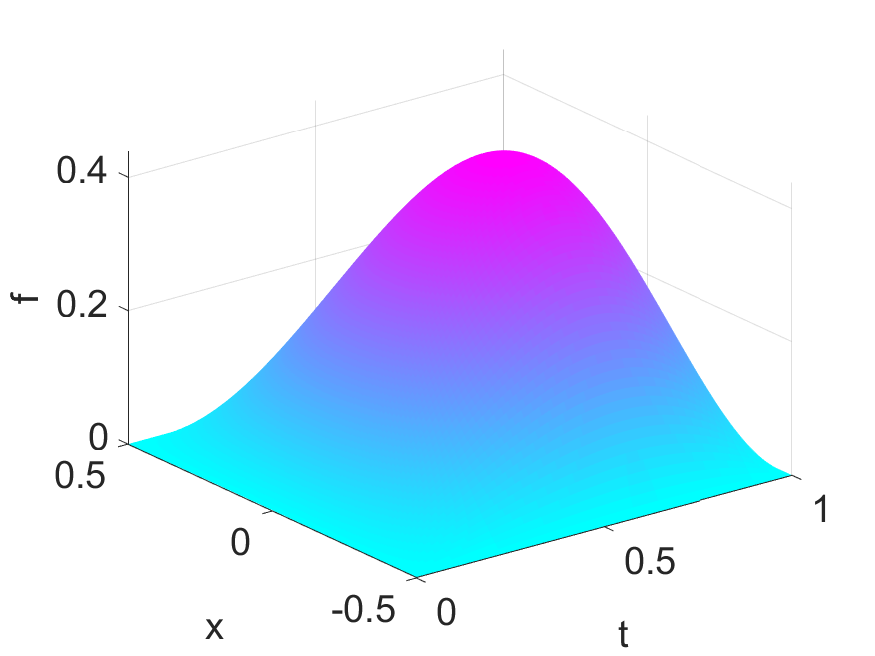} & \includegraphics[width=0.33\textwidth]{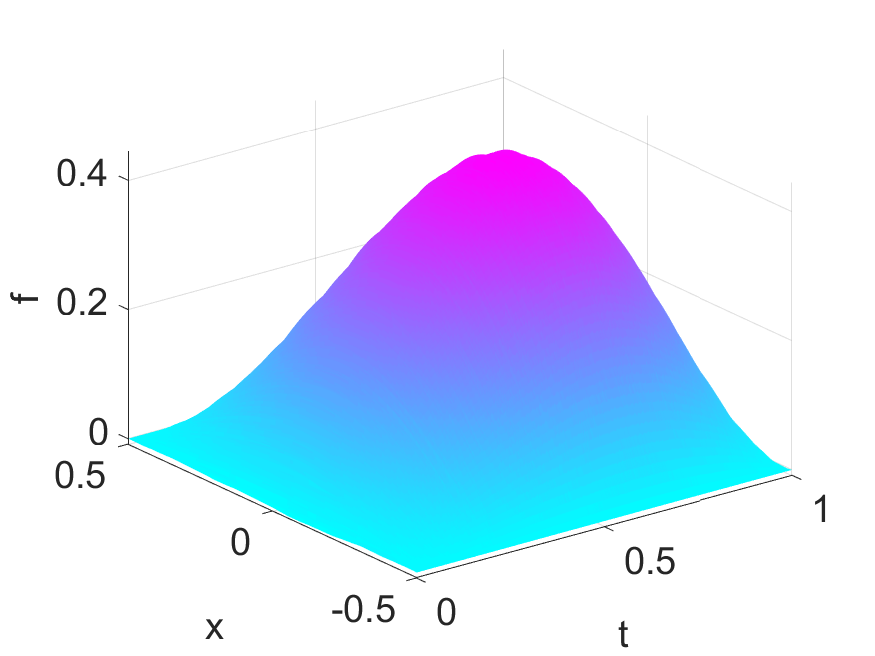} & \includegraphics[width=0.33\textwidth]{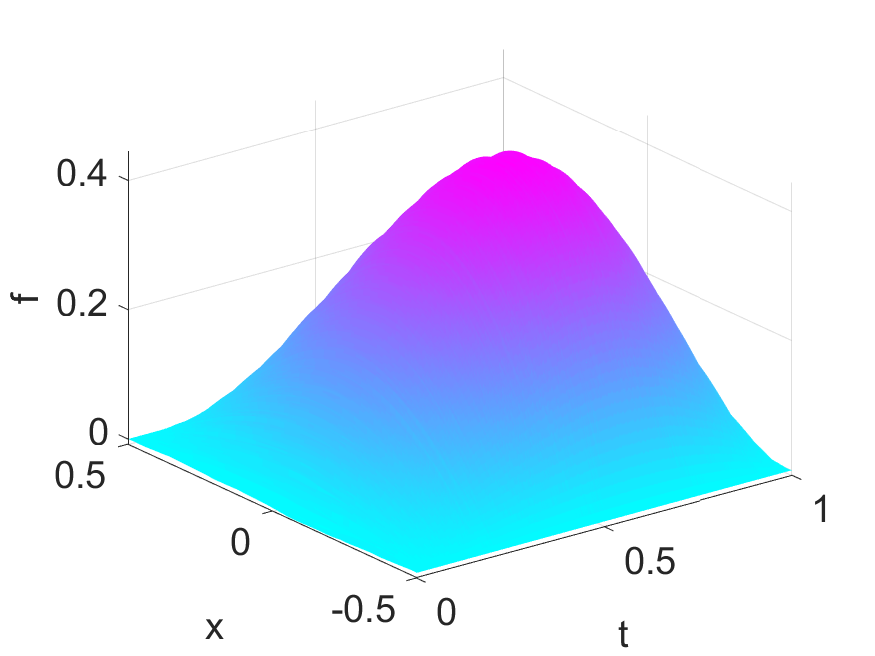}\\
       & \includegraphics[width=0.33\textwidth]{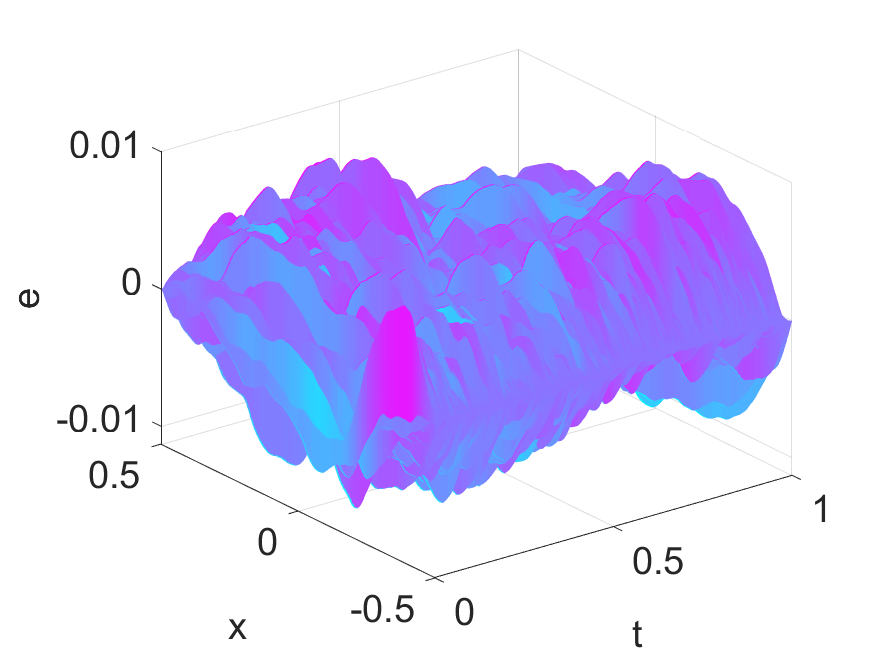} & \includegraphics[width=0.33\textwidth]{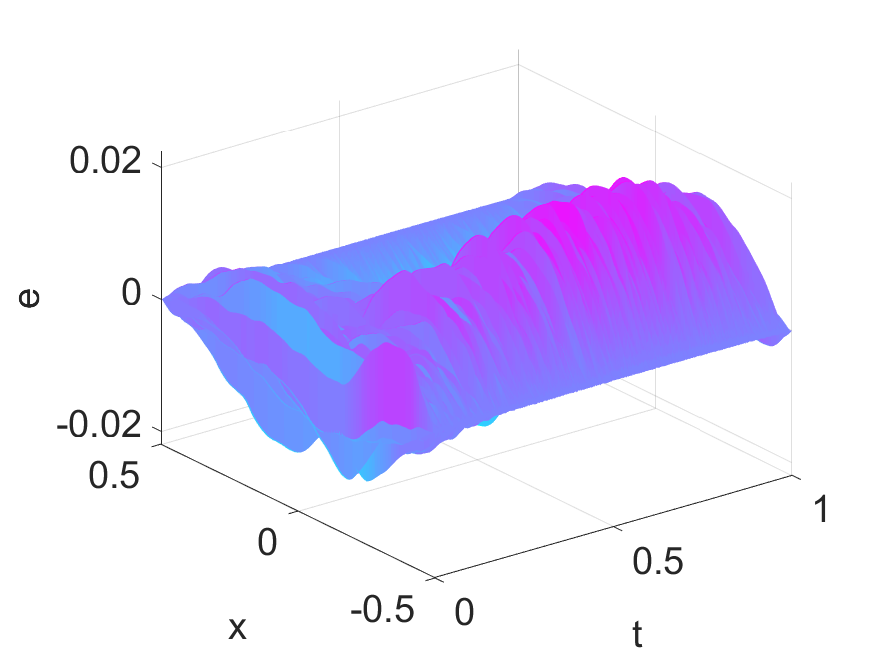}\\
       (a) exact & (b) $\varepsilon$=1e-2 & (c) $\varepsilon$=5e-2
  \end{tabular}
  \caption{Reconstructions and the pointwise errors for Example \ref{exam:dep}(i) with $\varepsilon$=1e-2 and $\varepsilon$=5e-2.\label{fig:dep1}}
\end{figure}

\begin{figure}[hbt!]
  \centering
  \setlength{\tabcolsep}{0pt}
  \begin{tabular}{cccc}
    \includegraphics[width=0.25\textwidth]{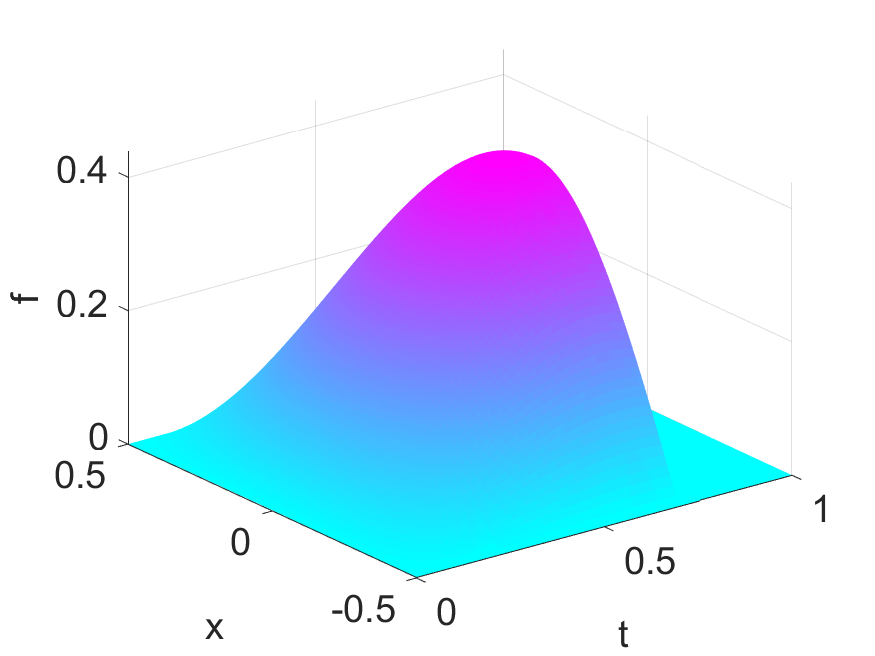} & \includegraphics[width=0.25\textwidth]{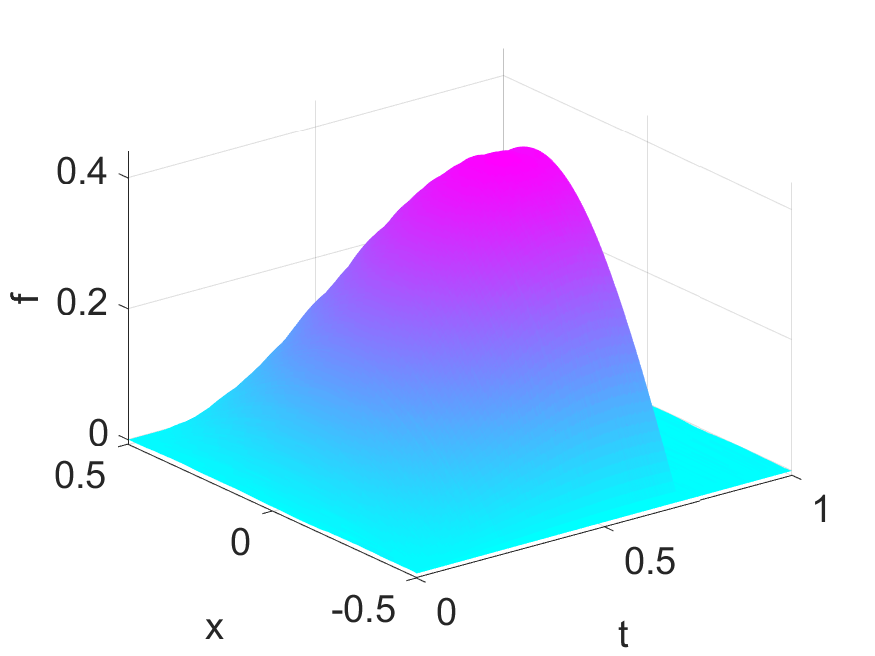} & \includegraphics[width=0.25\textwidth]{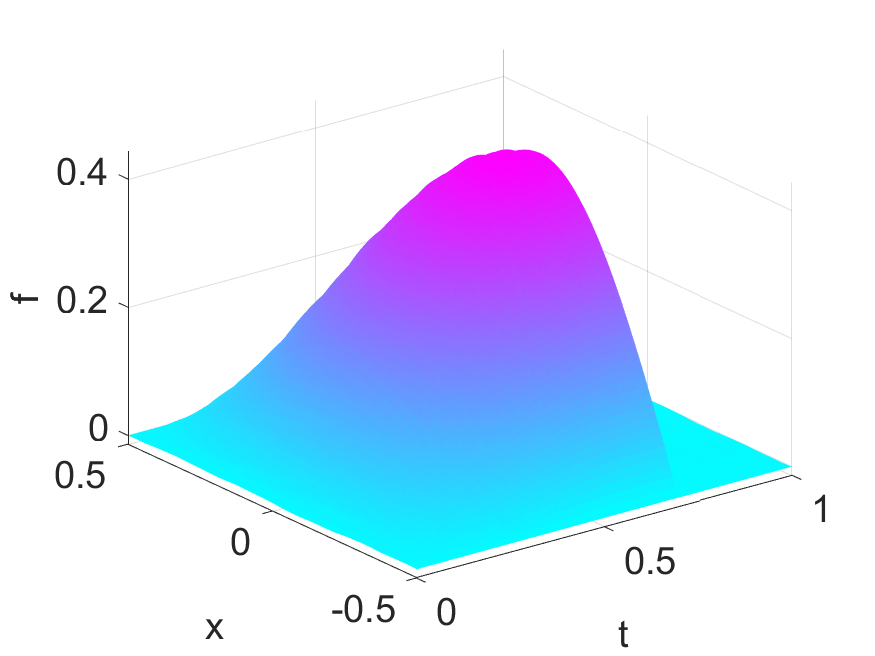}
     & \includegraphics[width=0.25\textwidth]{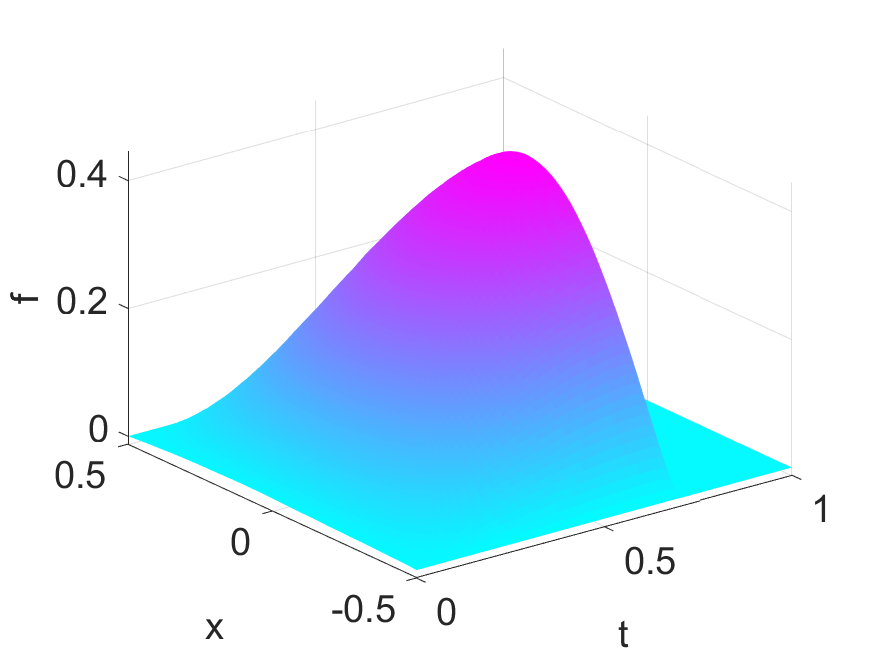}\\
       & \includegraphics[width=0.25\textwidth]{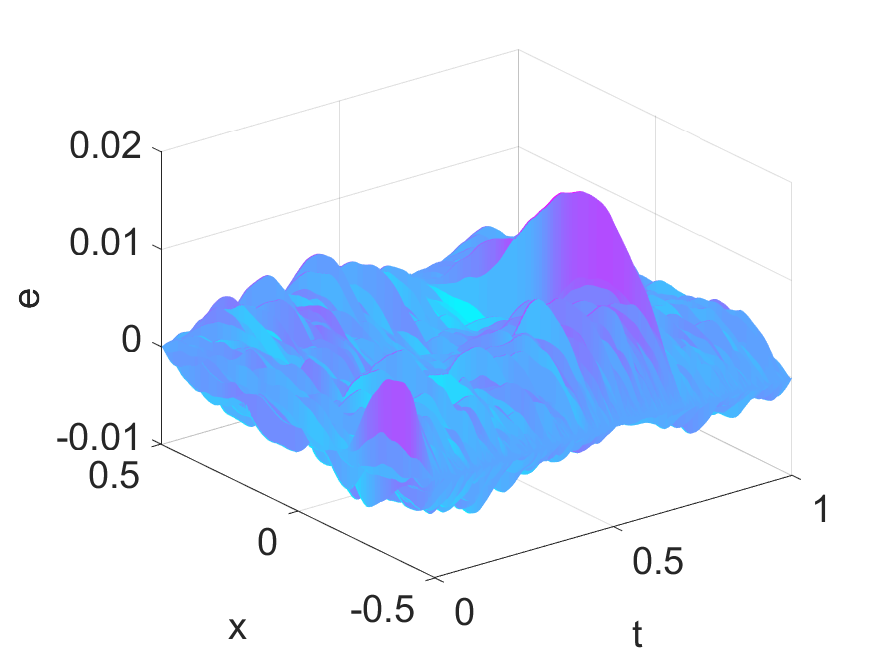} & \includegraphics[width=0.25\textwidth]{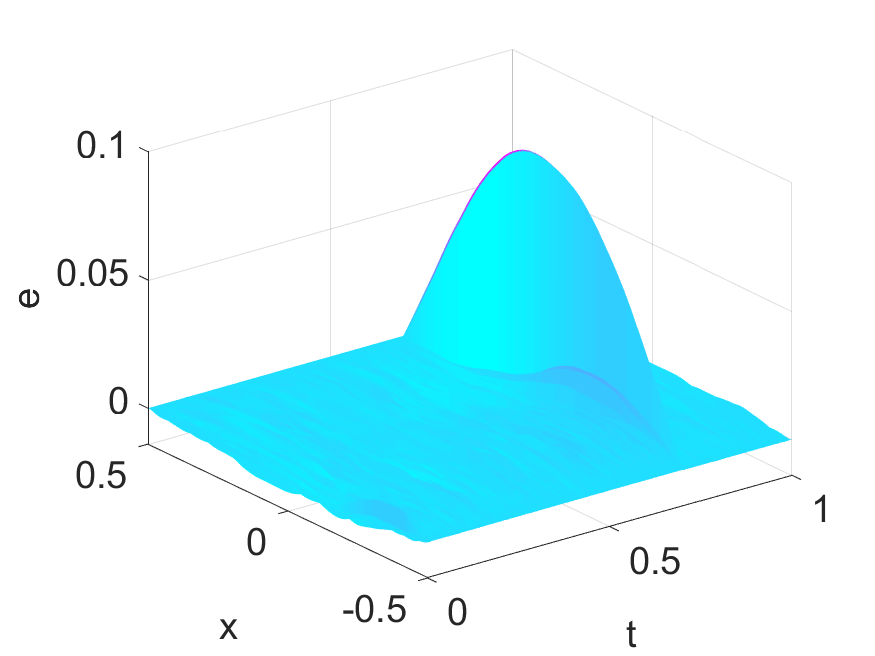}
     & \includegraphics[width=0.25\textwidth]{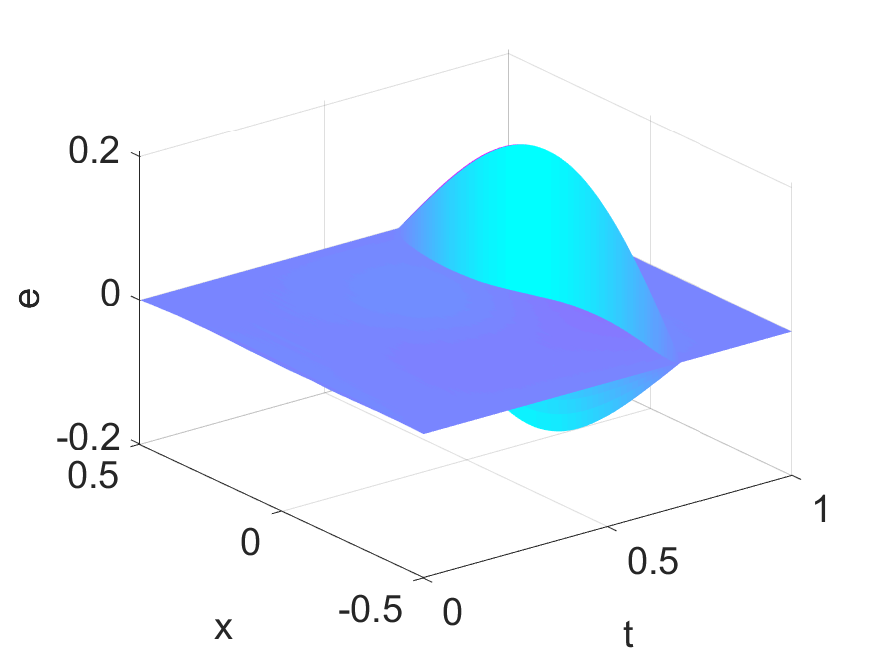}\\
       (a) exact & (b) $\alpha=0.25$ & (c) $\alpha=0.50$ & (d) $\alpha=0.75$
  \end{tabular}
  \caption{Reconstructions and the pointwise errors for Example \ref{exam:dep}(ii) with $\varepsilon$=1e-2.\label{fig:dep2}}
\end{figure}

\begin{figure}[hbt!]
\centering
\setlength{\tabcolsep}{0pt}
\begin{tabular}{cc}
  \includegraphics[width=0.48\textwidth]{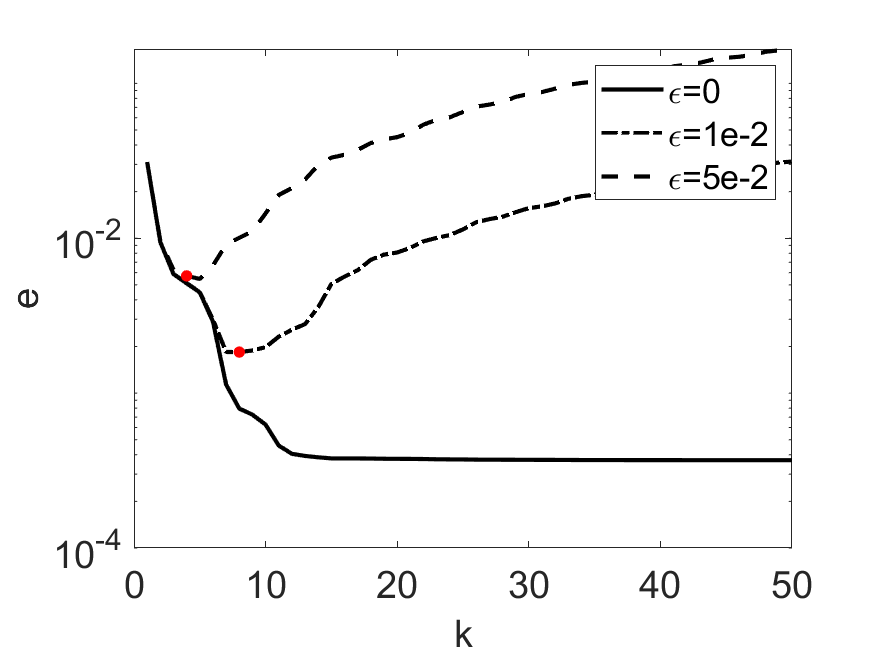} & \includegraphics[width=0.48\textwidth]{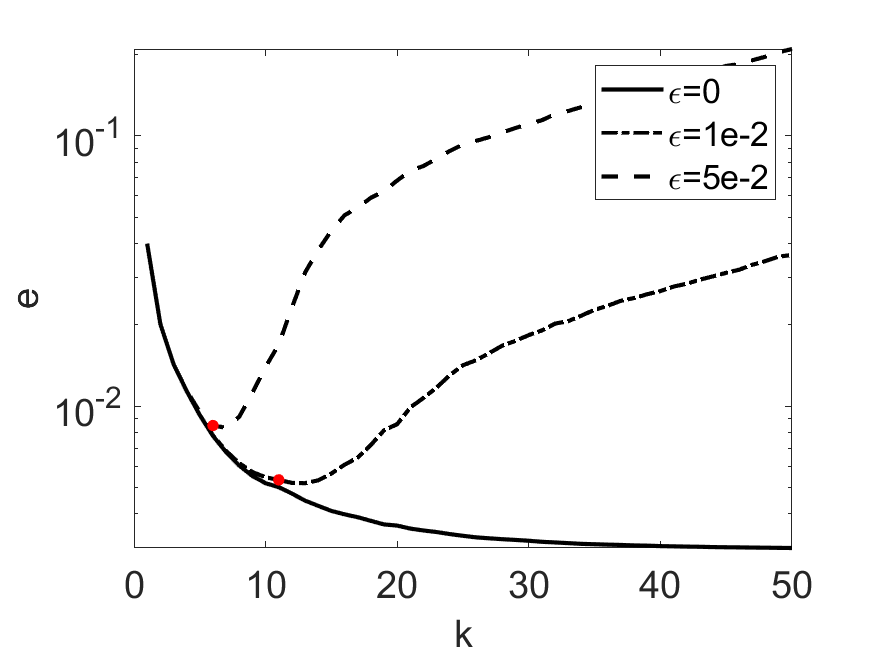}\\
  (a) $\alpha=0.25$ & (b) $\alpha=0.75$
\end{tabular}
\caption{The convergence of the error for Example \ref{exam:dep}(ii), where the red dots
indicate the stopping index determined by the discrepancy principle \eqref{eqn:dp}.\label{fig:conv}}
\end{figure}

In case (i), $f$ is smooth in time, but it is discontinuous for case (ii). The results for Example \ref{exam:dep}
are shown in Table \ref{tab:dep}. The results for case (i) are largely comparable with that for
Example \ref{exam:ind}, and all the observations remain valid; see also Fig. \ref{fig:dep1}.
The behavior of ISPn is largely independent of the fractional order $\alpha$,
due to the good regularity and compatibility of $f^\dag$. In sharp contrast, the results for case
(ii) exhibit a different trend: for a fixed noise level $\epsilon$, the reconstruction error $e$ increases
with the order $\alpha$, and also it takes more CG iterations to reach the convergence (see also Fig.
\ref{fig:conv}). This is attributed to the discontinuity in time of $f^\dag$
and the regularity of the adjoint $v$ in problem \eqref{eqn:fde-adj}: the temporal regularity of the adjoint $v$
increases steadily with $\alpha$, cf. Theorem \ref{thm:reg-adj}, which makes it increasingly
harder to approximate a discontinuous $f^\dag$. This is clearly visible from the error plots in Fig.
\ref{fig:dep2}, where the errors around the discontinuity dominate. This is especially
pronounced for $\alpha=0.50$ and $\alpha=0.75$.

\begin{table}[hbt!]
\centering
\caption{The reconstruction errors $e$ for Example \ref{exam:dep}.\label{tab:dep}}
\begin{tabular}{c|cccccc}
\toprule
case & $\alpha\backslash\varepsilon$ & 0 & 1e-3 & 5e-3 & 1e-2 & 5e-2\\
\midrule
   & 0.25 & 2.89e-5 (50) & 2.98e-4 (13) & 1.21e-3 (9) & 1.97e-3 (8) & 6.08e-3 (5)\\
(i)& 0.50 & 2.93e-5 (50) & 3.07e-4 (12) & 1.18e-3 (9) & 2.09e-3 (8) & 6.14e-3 (5)\\
   & 0.75 & 3.49e-5 (50) & 2.61e-4 (13) & 8.50e-4 (9) & 1.44e-3 (8) & 4.24e-3 (5)\\
\midrule
    & 0.25 & 3.69e-4 (50) & 4.51e-4 (13) & 1.12e-3 ( 9) & 1.84e-3 ( 8) & 5.71e-3 (4)\\
(ii)& 0.50 & 1.66e-3 (50) & 1.68e-3 (13) & 2.00e-3 (10) & 2.61e-3 ( 9) & 6.33e-3 (5)\\
    & 0.75 & 2.99e-3 (50) & 3.38e-3 (25) & 4.49e-3 (14) & 5.34e-3 (11) & 8.49e-3 (6)\\
\bottomrule
\end{tabular}
\end{table}

\subsubsection{Numerical results for ISPd}

Now we present two examples for ISPd, with the setting similar to Example \ref{exam:dep}.
\begin{example}\label{exam:dep-diri}
The diffusion coefficient $a$ is given by $a(x_1,x_2,t)=(1+\sin(\pi (x_1+\frac12))(\frac14-x_2^2))(1+\sin t)$, and consider two different source components.
\begin{itemize}
  \item[{\rm(i)}] $f^\dag (x_1,t) = \sin(x_1+\frac12)\pi t(T-t)e^t$.
  \item[{\rm(ii)}] $f^\dag (x_1,t) = \sin(x_1+\frac12)\pi t(T-t)e^t\chi_{[0,0.7]}(t)$.
\end{itemize}
\end{example}

\begin{table}[hbt!]
\centering
\caption{The reconstruction errors $e$ for Example \ref{exam:dep-diri}.\label{tab:dep-diri}}
\begin{tabular}{c|cccccc}
\toprule
case & $\alpha\backslash\varepsilon$ & 0 & 1e-3 & 5e-3 & 1e-2 & 5e-2\\
\midrule
   & 0.25 & 4.51e-3 (50) & 4.52e-3 (42) & 4.61e-3 (20) & 4.84e-3 (18) & 6.28e-3 (5)\\
(i)& 0.50 & 4.50e-3 (50) & 4.51e-3 (41) & 4.61e-3 (20) & 4.86e-3 (18) & 6.11e-3 (4)\\
   & 0.75 & 4.47e-3 (50) & 4.49e-3 (37) & 4.56e-3 (18) & 4.70e-3 (13) & 5.84e-3 (6)\\
\midrule
    & 0.25 & 3.87e-3 (50) & 3.88e-3 (36) & 3.97e-3 (20) & 4.20e-3 (17) & 5.47e-3 ( 4)\\
(ii)& 0.50 & 3.99e-3 (50) & 4.00e-3 (40) & 4.10e-3 (22) & 4.37e-3 (17) & 5.95e-3 ( 6)\\
    & 0.75 & 4.35e-3 (50) & 4.36e-3 (50) & 4.60e-3 (33) & 4.97e-3 (24) & 7.07e-3 (11)\\
\bottomrule
\end{tabular}
\end{table}

\begin{figure}[hbt!]
  \centering
  \setlength{\tabcolsep}{0pt}
  \begin{tabular}{cccc}
    \includegraphics[width=0.25\textwidth]{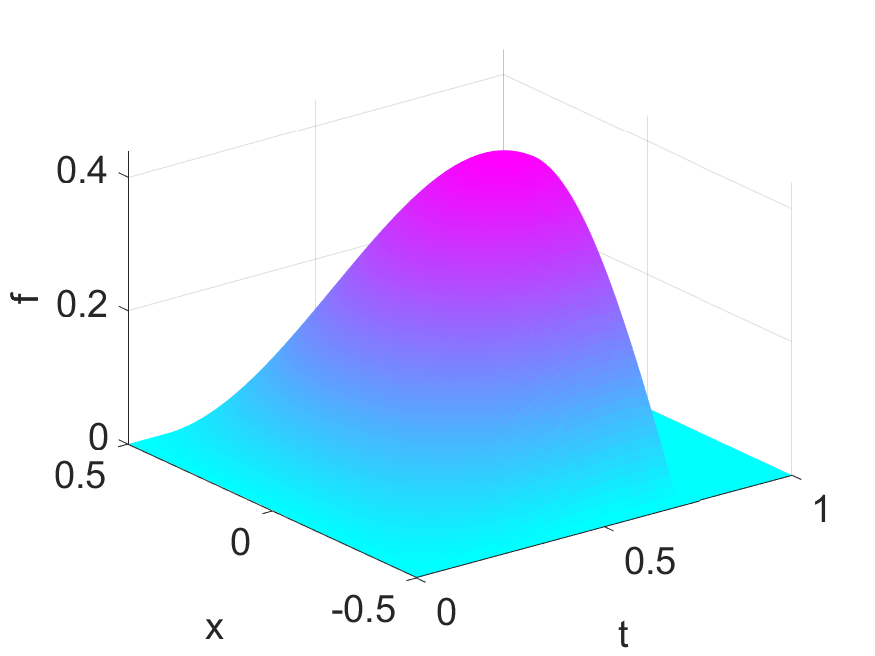} & \includegraphics[width=0.25\textwidth]{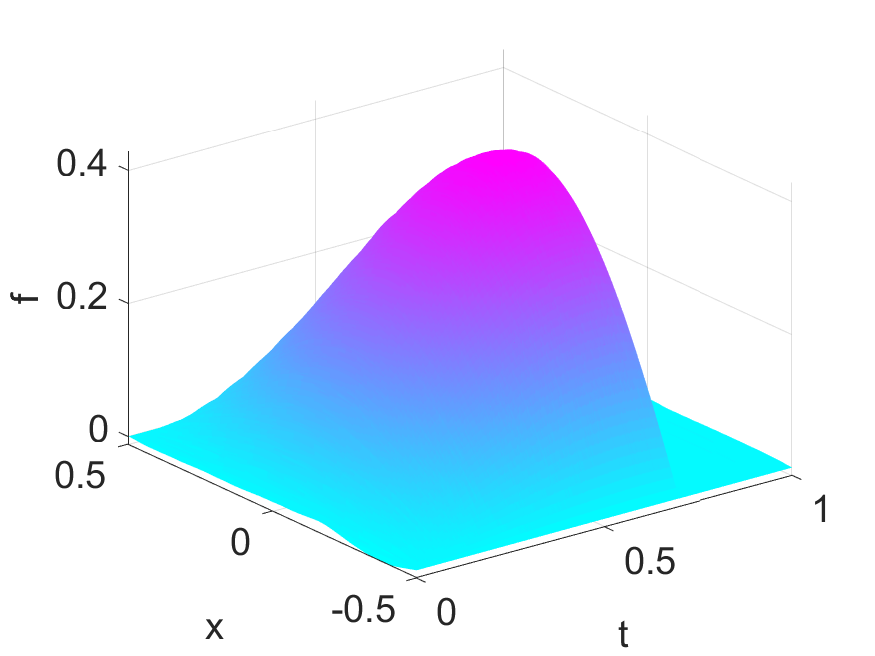} & \includegraphics[width=0.25\textwidth]{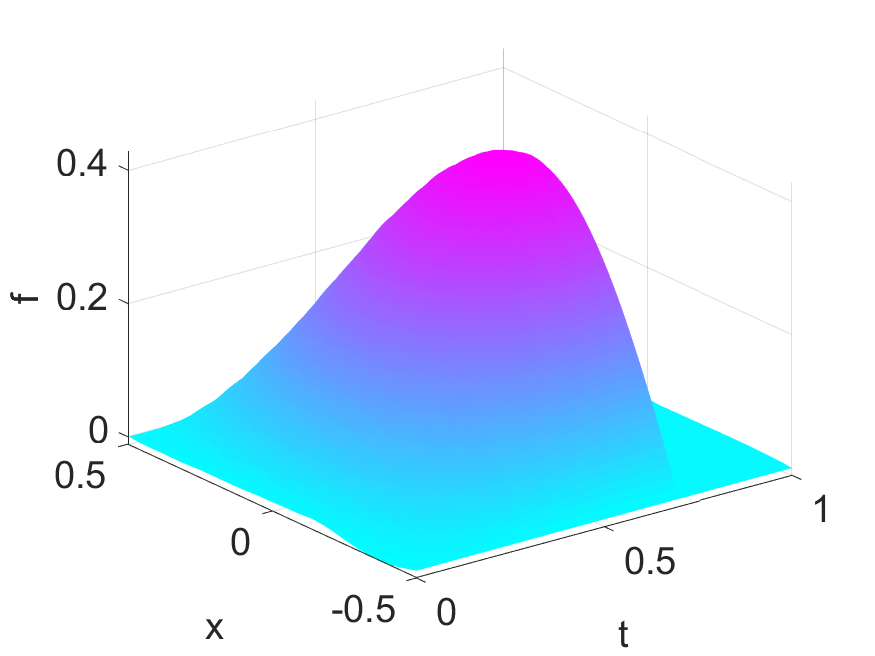}
     & \includegraphics[width=0.25\textwidth]{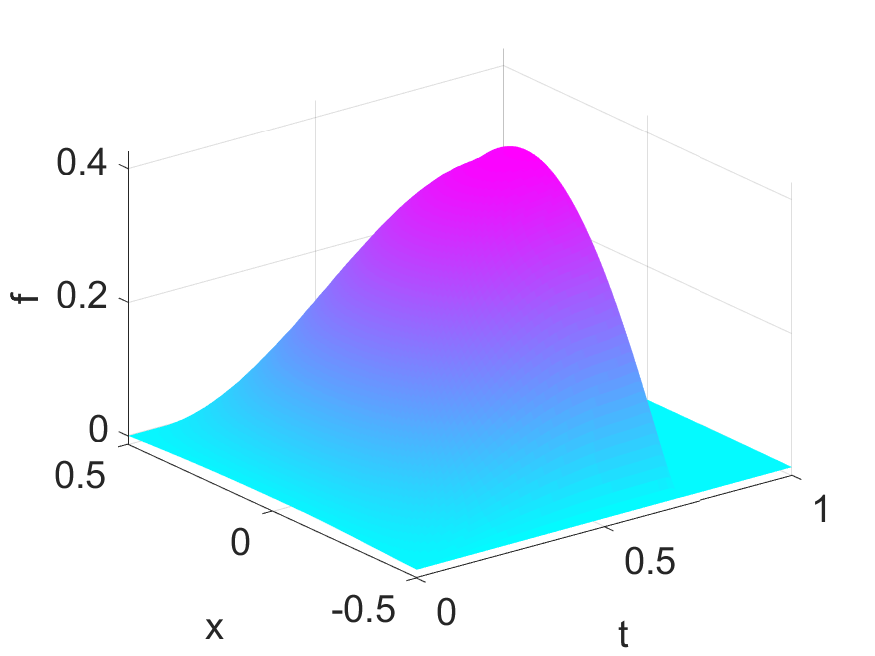}\\
       & \includegraphics[width=0.25\textwidth]{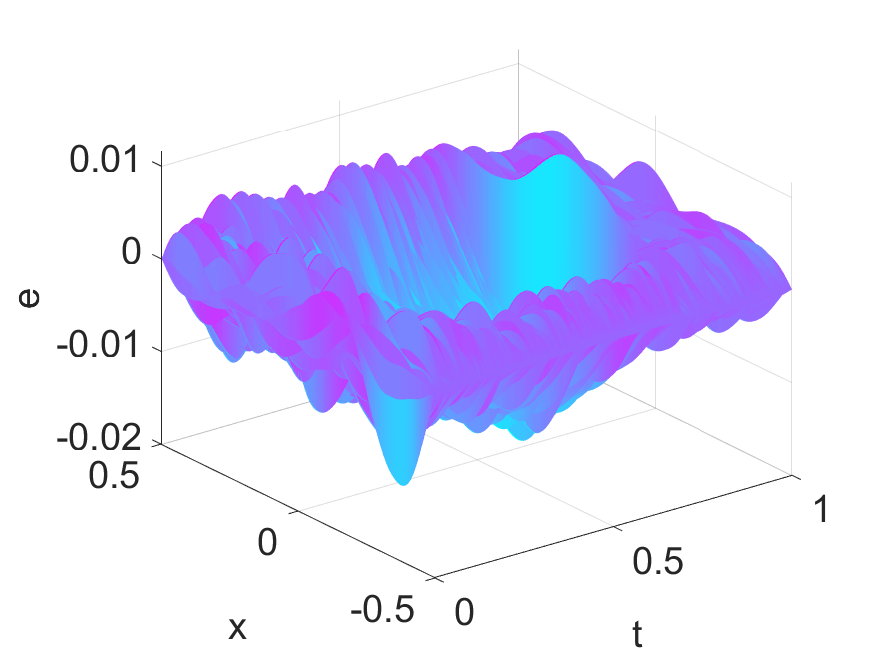} & \includegraphics[width=0.25\textwidth]{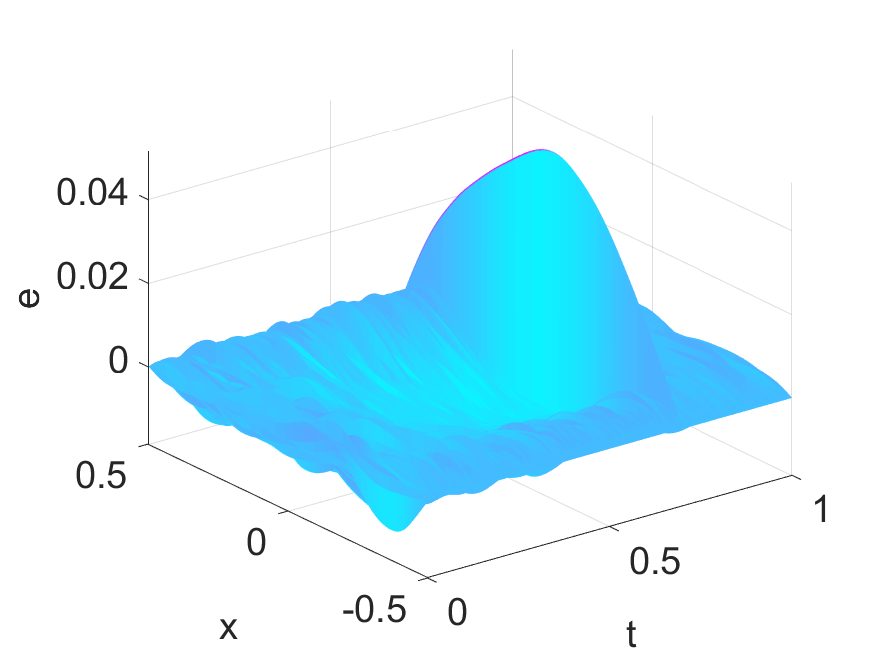}
     & \includegraphics[width=0.25\textwidth]{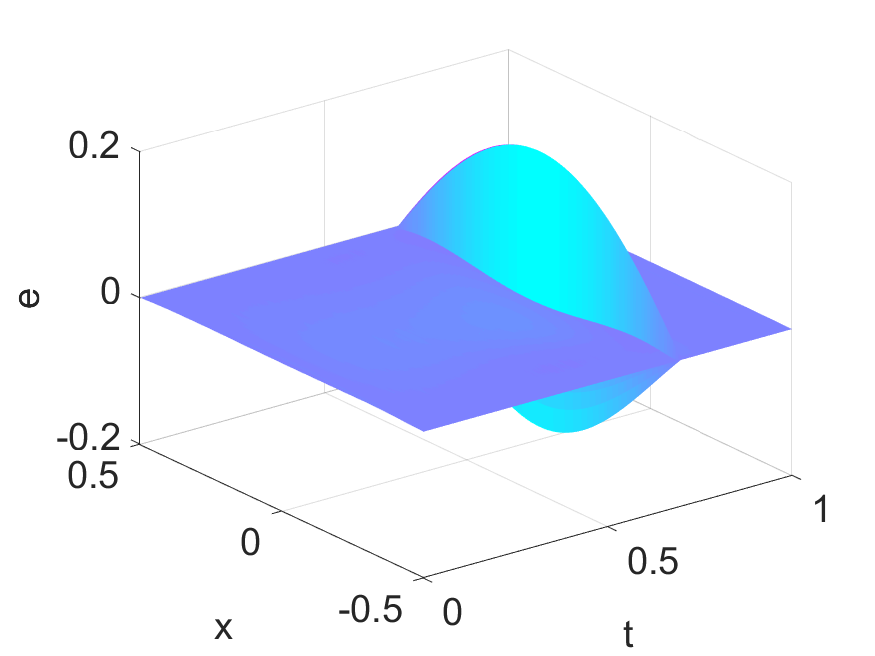}\\
       (a) exact & (b) $\alpha=0.25$ & (c) $\alpha=0.50$ & (d) $\alpha=0.75$
  \end{tabular}
  \caption{Reconstructions and the pointwise errors for Example \ref{exam:dep-diri}(ii) with $\varepsilon$=1e-2.\label{fig:dep2-diri}}
\end{figure}

Note that case (ii) does not satisfy the condition of Theorem \ref{thm:stability}. The numerical results
for Example \ref{exam:dep-diri} are shown in Table \ref{tab:dep-diri}, where the stopping index is taken
so that the reconstruction error $e$ is smallest (since the discrepancy principle \eqref{eqn:dp} does not
apply directly). The observations from Examples \ref{exam:ind} and \ref{exam:dep} are still valid, except
the algorithm takes more iterations to reach convergence. This might be due to the fact that the approximation
of the exact flux data (for the direct problem) is less accurate, which also limits the attainable accuracy
of the reconstruction for data with low noise level. The results for case (ii) show that for a fixed noise
level $\epsilon$, the error $e$ increases with $\alpha$, and also it takes more CG iterations to reach
convergence, due to the mismatch between the temporal regularity of $f$ and the gradient $J'(f)$. This is also clear from
the error plots in Fig. \ref{fig:dep2-diri}, where the errors around the discontinuity become increasingly
dominating as $\alpha$ increases.

\begin{figure}[hbt!]
  \centering
  \setlength{\tabcolsep}{0pt}
  \begin{tabular}{ccc}
    \includegraphics[width=0.33\textwidth]{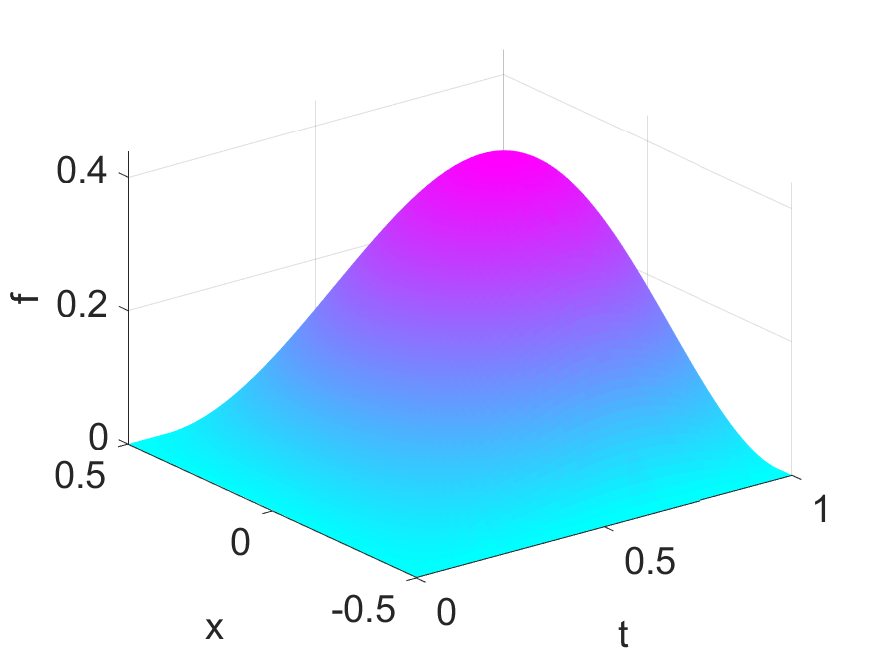} & \includegraphics[width=0.33\textwidth]{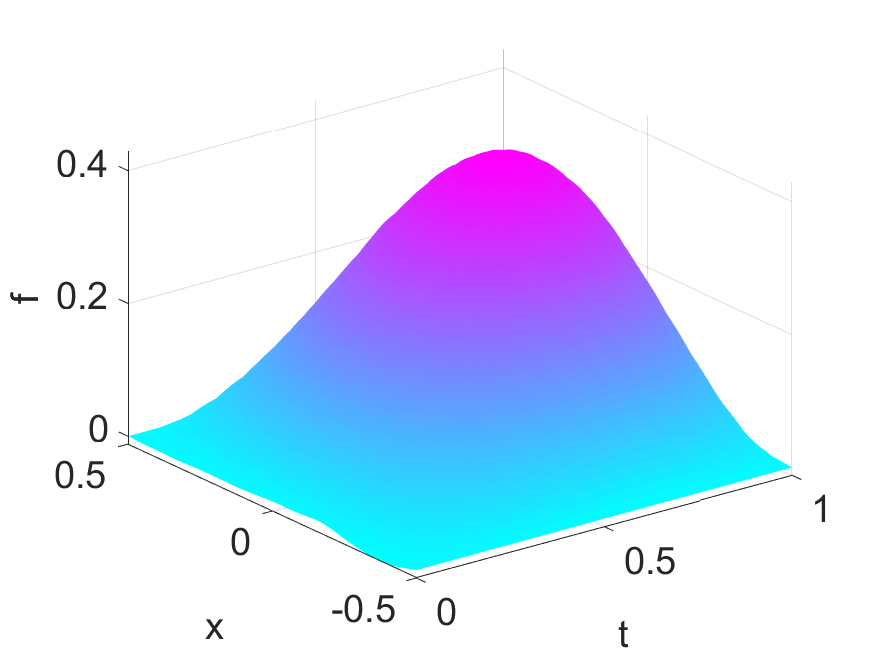} & \includegraphics[width=0.33\textwidth]{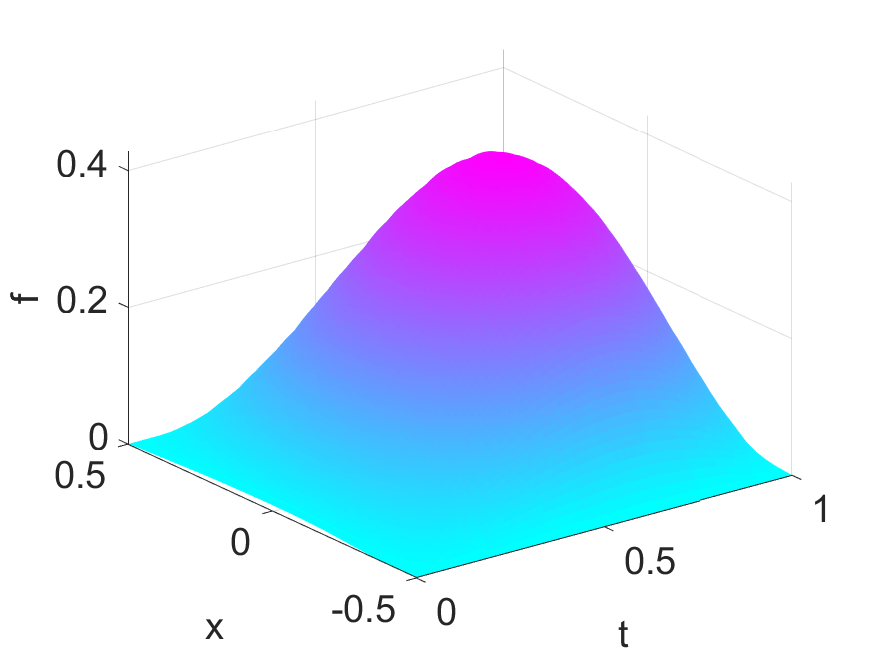}\\
       & \includegraphics[width=0.33\textwidth]{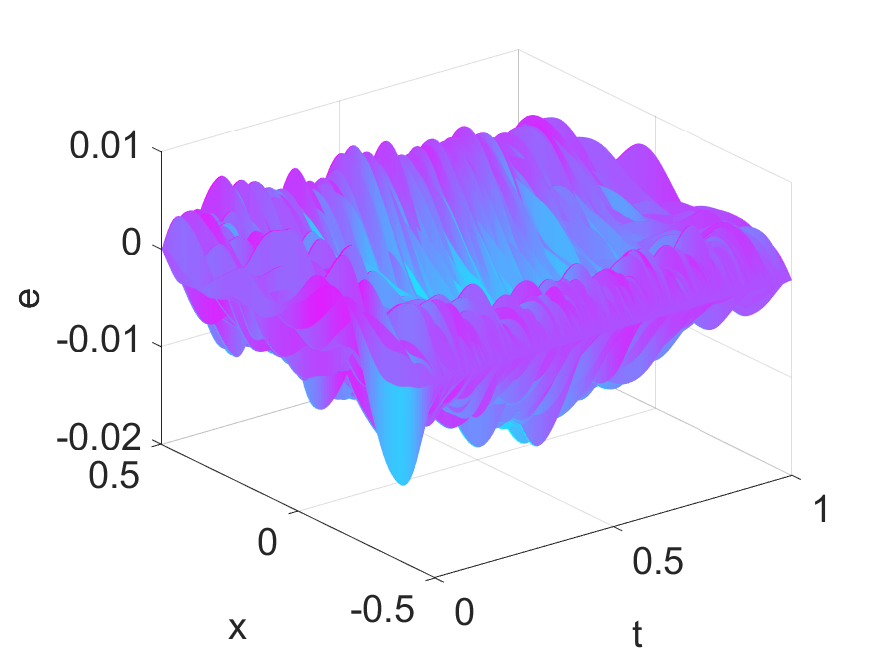} & \includegraphics[width=0.33\textwidth]{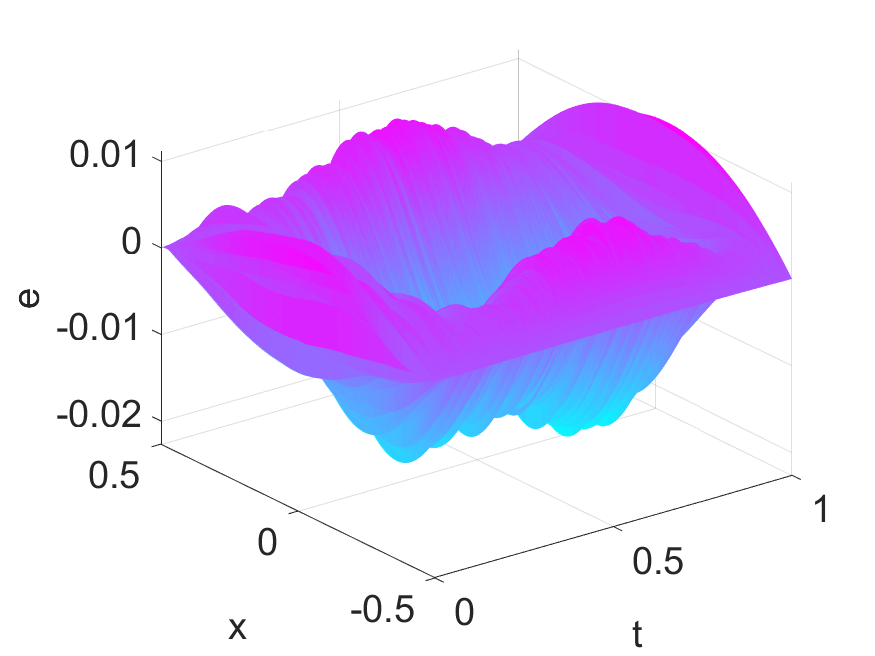}\\
       (a) exact & (b) $\varepsilon$=1e-2 & (c) $\varepsilon$=5e-2
  \end{tabular}
  \caption{Reconstructions and the pointwise errors for Example \ref{exam:dep-diri}(i) with $\varepsilon$=1e-2 and $\varepsilon$=5e-2.\label{fig:dep1-diri}}
\end{figure}

These numerical results indicate that indeed it is feasible to recover
a space-time dependent source from the lateral boundary observation in a cylindrical
domain for both time-independent and dependent diffusion coefficients, and standard
regularization techniques, e.g., conjugate gradient method (when equipped with the discrepancy
principle \eqref{eqn:dp}), can deliver accurate reconstructions for both exact and noisy data. This provides
numerical evidences to the theoretical results in Theorems \ref{thm:uniqueness} and \ref{thm:stability}.

\bibliographystyle{abbrv}
\bibliography{frac}
\end{document}